\documentclass[]{article}
\usepackage{graphicx,amsfonts,amsmath,amsthm,amssymb,amscd,longtable}

\def\tr{{\rm tr}}
\def\PSL{{\rm PSL}(2,{\Bbb C})}
\def\P{{\cal P}}

\newtheorem{theorem}{Theorem}
\newtheorem{lemma}{Lemma}
\newtheorem{proposition}{Proposition}

\title{All discrete ${\cal RP}$ groups whose generators have real
traces}
\author
{Elena Klimenko\footnote{Partially supported by the NSF grant DMS-9971660.}
\and
Natalia Kopteva}

\date{}

\begin{document}

\maketitle

\begin{abstract}
In this paper we give necessary and sufficient conditions
for discreteness of a subgroup of ${\rm PSL}(2,{\Bbb C})$ generated by
a hyperbolic element and
an elliptic one of odd order with non-orthogonally intersecting axes.
Thus we completely determine two-generator non-elementary Kleinian
groups without invariant plane with real traces of the generators
and their commutator. We also give a list of all parameters that
correspond to such groups.
An interesting corollary
of the result is that the group of the minimal known volume
hyperbolic orbifold ${\Bbb H}^3/\Gamma_{353}$ has real parameters.
\end{abstract}


\noindent
Secondary: 20H10, 22E40, 57M60, 57S30.

\normalsize

\section{Introduction}

The group of all M\"obius transformations of the extended complex
plane $\overline {\Bbb C}={\Bbb C}\cup\lbrace\infty\rbrace$ is isomorphic to
${\rm PSL}(2,{\Bbb C})={\rm SL}(2,{\Bbb C})/\lbrace\pm I\rbrace$.
The Poincar\'e extension gives the action of this group
(as the group of all orientation preserving isometries)
on hyperbolic 3-space
$$
{\Bbb H}^3=\lbrace(z,t)\ | \ z\in{\Bbb C},\ t>0\rbrace
$$
with the Poincar\'e metric
$$
\displaystyle ds^2=\frac{|dz|^2+dt^2}{t^2}.
$$

We study the class of $\cal RP$ {\it groups} (two-generator
groups with real parameters):
$$
{\cal RP}=\lbrace\Gamma=\langle f,g\rangle \ | \
f,g\in{\rm PSL}(2,{\Bbb C});\ \beta,\beta',\gamma\in {\Bbb R}\rbrace,
$$
where $\beta=\beta(f)={\rm tr}^2f-4$, $\beta'=\beta(g)={\rm tr}^2g-4$,
$\gamma=\gamma(f,g)={\rm tr}[f,g]-2$.

Recall that an element $f\in\PSL$ with real $\beta(f)$
is
{\it elliptic}, {\it parabolic}, {\it hyperbolic}, or {\it $\pi$-loxodromic}
according to whether
$\beta(f)\in[-4,0)$, $\beta(f)=0$, $\beta(f)\in(0,+\infty)$, or
$\beta(f)\in(-\infty,-4)$.
If $\beta(f)\notin[-4,\infty)$,
then $f$ is called {\it strictly loxodromic}.
Among all
strictly loxodromic elements only $\pi$-loxodromics have
real $\beta(f)$.

\medskip\noindent
{\bf 1.1. Elementary $\cal RP$ groups.}
Following \cite{B}, we call a subgroup $\Gamma$ of
${\rm PSL}(2,{\Bbb C})$ {\it elementary} if there
exists a finite $\Gamma$-orbit in ${\Bbb H}^3\cup\overline{\Bbb C}$.
Elementary groups are studied in~\cite{B,M}. Among the $\cal RP$
groups, the following are elementary:
\begin{itemize}
\item[(1)] Both generators are elliptic with intersecting axes.
\item[(2)] Both generators are elliptic of order~$2$ with mutually
orthogonal skew axes or with disjoint axes lying in a hyperbolic plane.
\item[(3)] One generator is hyperbolic or $\pi$-loxodromic and the other
is elliptic of order $2$ whose axes intersect orthogonally.
\item[(4)] The generators (with real $\beta$, $\beta'$, and $\gamma$)
share a fixed point in $\overline{\Bbb C}$.
\end{itemize}
The list of all {\it discrete} elementary groups (not necessarily from
the $\cal RP$ class) is given in \cite{B}, so there is
no need to consider discrete elementary ${\cal RP}$ groups here.

\medskip\noindent
{\bf 1.2. $\cal RP$ groups with invariant hyperbolic plane.}
\begin{itemize}
\item[(1)] Any two-generator Fuchsian group, i.e., conjugate to a
subgroup of ${\rm PSL}(2,{\Bbb R})$, is an $\cal RP$ group. Discreteness
conditions for such groups have been well-studied (see~\cite{Gi,GiM,R},
for a complete list of references see~\cite{FR}).
\item[(2)] Any non-Fuchsian two-generator subgroup of
${\rm PSL}(2,{\Bbb C})$ with invariant hyperbolic plane is also an
$\cal RP$ group. Such a group contains elements that reverse orientation
of the invariant plane and interchange the half-spaces bounded by it.
Discreteness criteria and a classification theorem for all such
discrete groups are given in \cite{KS}.
\end{itemize}
In \cite{GGM}, the space of discrete $\cal RP$ groups with parameters
$(\beta,-4,\gamma)$ is analyzed. Since one of the generators is elliptic
of order $2$ and parameters $\beta$ and $\gamma$ are real, it follows from
\cite[Theorems 1(i) and 2(i)]{KK} that all such groups are either elementary
or have an invariant plane.

Our paper is devoted to the following class:

\medskip\noindent
{\bf 1.3. Truly spatial $\cal RP$ groups.}
We call a subgroup of ${\rm PSL}(2,{\Bbb C})$ {\it truly spatial}, if it
is not elementary and has no invariant hyperbolic plane.
All truly spatial ${\cal RP}$ groups are characterized by the following theorem
(see also Table~1 in~\cite{KK}):
\begin{theorem}\label{th2}
{\rm {\cite[Theorem 4]{KK}}}
Let \ $\Gamma=\langle f,g\rangle$ be an
${\cal RP}$ group.
$\Gamma$ is a truly spatial group if and only if
\begin{center}
$(-1)^k\gamma<(-1)^{k+1}\beta\beta'/4$ with
$\gamma\not=0$, $\beta\not=-4$, and $\beta'\not=-4$,
\end{center}
where $k\in\lbrace 0,1,2 \rbrace$ is the number of $\pi$-loxodromic
elements among $f$ and $g$.
\end{theorem}

The current paper  is the last in a series of earlier papers
\cite{K2}--\cite{KK05} that gathered together give
criteria (necessary and sufficient conditions) for discreteness of
all truly spatial $\cal RP$ groups with real traces of
both generators. Thus, we complete
the description of all discrete $\cal RP$ groups with
non-$\pi$-loxodromic generators. The result in terms of parameters
can be found in Table~\ref{table_param} in Appendix, where all real triples
$(\beta,\beta',\gamma)$ for non-elementary discrete $\cal RP$
groups without invariant plane with non-$\pi$-loxodromic generators
are listed.

\medskip
Our main result is Theorem~A below that gives criteria for
discreteness of truly spatial $\cal RP$ groups generated by a hyperbolic
element and an elliptic one of odd order with intersecting
axes.
The discrete groups from this theorem correspond to
rows 21--41 of Table~\ref{table_param} in Appendix.
In the theorem we assume without loss of generality
that the elliptic generator is primitive; if not, we can replace it by
its primitive power of the same order (cf. \cite{KK}).

\medskip
\noindent
{\bf Theorem A.}\label{thm1}
{\it
Let $f$ be a primitive elliptic element of odd order $n$, $n\geq 3$,
$g$ be a~hyperbolic element, and let their axes intersect
non-orthogonally. Then
\begin{itemize}
\item[$(1)$] There exist elements $h_1,h_2\in \PSL$ such that
$h_1^2=gfg^{-1}f$,
$h_2^2=f^{(n-1)/2}g^{-1}f^{-1}gf^{-(n+1)/2}gf^{-1}g^{-1}$,
$(h_1f^{-1})^2=1$, and $h_2gfg^{-1}$ is an elliptic element whose
axis intersects the axis of $f$. There exists also an element
$h_3\in\PSL$ such that
$h_3^2=f^{(n-1)/2}g^{-1}h_1^{-1}gf^{-(n-3)/2}h_1^{-1}$ and
$h_3h_1$ is an elliptic element whose axis intersects the axis of
$f$. \item[$(2)$] $\Gamma=\langle f,g\rangle$ is discrete if and
only if one of the following conditions holds:
\begin{itemize}
\item[$(i)$] $h_1$ is hyperbolic, parabolic, or a primitive
elliptic element of even order $2m$ $(2/n+1/m<1)$ and $h_2$ is
hyperbolic, parabolic, or a~primitive elliptic element of even
order $2l$, $l\geq 2$;
\item[$(ii)$] $h_1$ and $h_2$ are rotations
through angles of $\pi/m$ and $\pi/l$, respectively, and the ordered
triple $\{n,m,l\}$ is one of the following: $\{5,2,5/2\}$;
$\{3,m,m/3\}$, $m$ is an integer, $m\geq 4$, $(m,3)=1$; $\{n,3,n/3\}$,
$n\geq 5$, $(n,3)=1$; $\{3,5,3/2\}$; $\{5,2,5/3\}$; $\{3,5,5/4\}$;
$\{5,3,5/4\}$; $\{5,2,3/2\}$;
\item[$(iii)$] $n=3$, $h_1$ is
hyperbolic, $h_2$ is a rotation through an angle of $4\pi/k$,
$h_2gfg^{-1}$ is a primitive elliptic element of order $p$ and
the ordered pair
$\{k,p\}$ is one of the following: $\{6,5\}$; $\{k,3\}$, $k$ is an
integer, $k\geq 7$, $(k,4)\leq 2$;
\item[$(iv)$] $h_1$ is
hyperbolic and $h_2$ is the cube of a primitive elliptic element
$\tilde h_2$ of order $2n$
so that  $\tilde h_2^2gfg^{-1}$ is a rotation through
angle $4\pi/k$, where the ordered pair $\{n,k\}$ is one of the
following:
$\{5,5\}$; $\{n,4\}$, $n\geq 5$, $(n,3)=1$;
\item[$(v)$] $h_1$ is a primitive elliptic
element of odd order $\widetilde m$ $(1/n+1/\widetilde m<1/2)$
and $h_3$ is
hyperbolic, parabolic, or a primitive elliptic element of even
order $2k$, $k\geq 2$;
\item[$(vi)$] $n=3$, $h_1$ and $h_3$ are
primitive elliptic elements of the same order $\widetilde m$,
$\widetilde m\geq 7$, $(\widetilde m,2)=1$;
\item[$(vii)$] $h_1$ is the square of a primitive elliptic element
$\tilde h_1$ of order $n$, $n\geq 7$, and $h_4$ is hyperbolic,
parabolic, or a primitive elliptic element of even order $\geq 4$,
where $h_4$ is defined
as follows: $h_4^2=f^{(n-3)/2}g^{-1}\tilde h_1 gf^{-(n+1)/2}\tilde
h_1 f^{-1}$ and $h_4f\tilde h_1^{-1}$ is an elliptic element whose
axis intersects the axis of $f$.
\end{itemize}
\end{itemize}
}

\medskip
The content of the paper is mainly the proof of this theorem
that is given in Sections~3--6.
In Section~8 we reformulate Theorem~A in terms of parameters
to get Theorem~B.

\medskip\noindent
{\bf Minimal volume hyperbolic $3$-orbifold.}
Let $\Gamma_{353}$ be a ${\Bbb Z}_2$-extension of the orientation
preserving index $2$ subgroup of the group generated by reflections
in the faces of the hyperbolic tetrahedron
with Coxeter diagram 3--5--3. This group is known to have minimal
covolume (=0.03905...) among all arithmetic groups \cite{CF} and
among all groups with a tetrahedral subgroup (the symmetry group
of a regular tetrahedron)
and groups containing elliptic elements of order $p\geq 4$
\cite{GM1}. The latter fact is a great evidence of the Gehring-Martin
conjecture that $\Gamma_{353}$ has minimal covolume among all
Kleinian groups. In Section~7 we prove that $\Gamma_{353}$ is an
${\cal RP}$ group.

\bigskip
\noindent
{\bf Acknowledgements.}
The authors would like to thank the Edinburgh Mathematical
Society for its financial support which enabled the first
author to visit Heriot-Watt University, and personally
Prof.~J.~Howie for his splendid hospitality during the
preparation of the paper and many interesting discussions.

We also thank the referee who gave us an idea of how to
improve the introduction to this paper.

\section{Preliminaries}

For the benefit of the reader, we will list some preliminary
definitions and results that we will need in this paper.

For basic definitions in discrete groups we refer the
reader to~\cite{M}. All formulas of hyperbolic geometry that we
use in calculations were taken from~\cite{Fen}.

We also use the following powerful instruments in our proof:
discrete extensions of tetrahedral groups and minimal distances
between axes of elliptic elements in discrete groups.

Denote by $G_T$ the group generated by reflections in the faces of a
compact Coxeter tetrahedron $T$. All discrete extensions of its
orientation preserving index 2 subgroup $\Delta_T$ were obtained
by Derevnin and Mednykh in~\cite{DM}.

The minimal distance $\rho_{min}(p,q)$
between the axes of two elliptic elements
of orders $p$ and $q$ in a discrete group
were obtained by Gehring and Martin in~\cite{GM}
and  Gehring, Marshall, and
Martin in~\cite{GMM}. Table~\ref{table_min} shows
$\cosh\rho_{min}(p,q)$ for $p,q\leq 7$.

\begin{table}[htbp]
\begin{center}
\caption{$\cosh\rho_{min}(p,q)$, $p,q\leq 7$}\label{table_min}
\begin{tabular}{|c|c|c|c|c|c|c|c|}
\hline
 & 2 & 3& 4 & 5 & 6 & 7  \\
\hline
2 & 1.000 & 1.019 & 1.088 & 1.106 & 1.225 & 1.152  \\
3 & 1.019 & 1.079 & 1.155 & 1.376 & 1.155 & 1.198  \\
4 & 1.088 & 1.155 & 1.366 & 1.203 & 1.414 & 1.630  \\
5 & 1.106 & 1.376 & 1.203 & 1.447 & 1.701 & 1.961  \\
6 & 1.225 & 1.155 & 1.414 & 1.701 & 2.000 & 2.305  \\
7 & 1.152 & 1.198 & 1.630 & 1.961 & 2.305 & 1.656  \\
\hline
\end{tabular}
\end{center}
\end{table}

In fact, for $p\geq 7$ and $q\geq 2$ the minimal distance
is determined by the following formula~\cite{GMM}:
\begin{equation}\label{mindist}
2\sin(\pi/p)\sin(\pi/q)\cosh \rho_{min}=\left\{
\begin{array}{lll}
1 & {\rm if} & q\neq 3, q\neq p,\\
\cos(\pi/p) & {\rm if} & q=3,\\
\cos(2\pi/p) & {\rm if} & q=p.\\
\end{array}
\right.
\end{equation}

We will also need the following geometric Lemmas.

\begin{figure}[htbp]
\centering
\begin{tabular}{cc}
\includegraphics[width=4 cm]{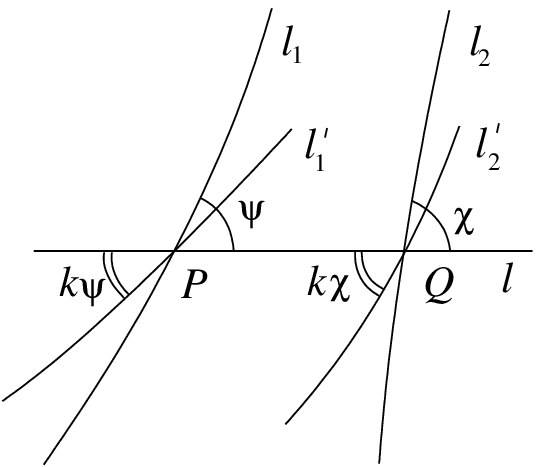}\qquad\quad &
\quad\qquad \includegraphics[width=3.8 cm]{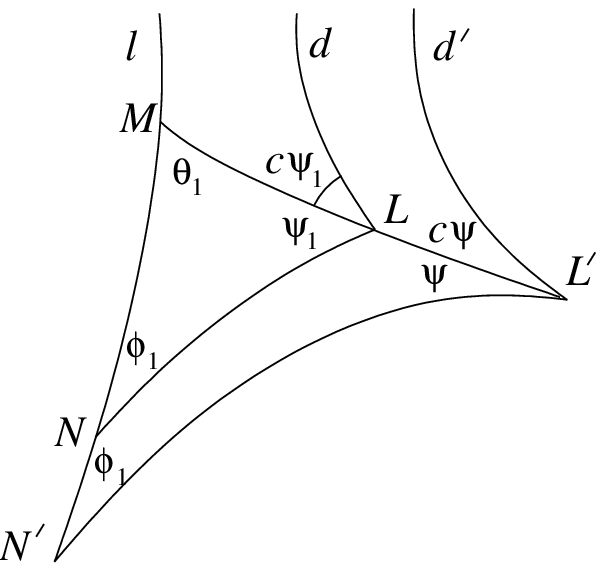}\\
(a)\qquad\quad & \qquad (b)\\
\end{tabular}
\caption{}\label{l1l2}
\end{figure}

\begin{lemma}\label{lemma_lines}
Let $l_1$ and $l_2$ be two disjoint lines
in ${\Bbb H}^2$ and let $l$ be a transversal of them meeting them in
points $P$ and $Q$ with corresponding angles of $\psi$ and $\chi$,
respectively, where $\psi<\chi<\pi/2$.
Let $l'_1$ and $l'_2$ be new lines passing through $P$ and
$Q$ so that their corresponding angles are of $k\psi$ and $k\chi$,
respectively,
$0<k<1$ (see Figure~\ref{l1l2}a). Then $l'_1$ and $l'_2$ are also
disjoint.
\end{lemma}

\noindent
{\it Proof.}
By contradiction. Suppose $l'_1$ and $l'_2$ are parallel or
intersecting. Then
\begin{equation}\label{leq1}
\cosh PQ\leq\frac{1-\cos k\psi\cos k\chi}{\sin k\psi\sin k\chi}
\end{equation}
Note that (\ref{leq1}) is an identity for parallel lines.
Now prove that
$$
f(k)=\frac{1-\cos k\psi\cos k\chi}{\sin k\psi\sin k\chi},
\qquad 0<\psi<\chi<\pi/2,
$$
is increasing on $(0,1]$.
Indeed,
\begin{eqnarray*}
f'(k)(\sin k\psi\sin k\chi)^2
&=& (\psi\sin k\psi\cos k\chi+\chi\sin k\chi\cos k\psi)
\sin k\psi\sin k\chi+\\
&&(-1+\cos k\psi\cos k\chi)
(\psi\cos k\psi\sin k \chi+\chi\cos k\chi\sin k\psi)\\
&=& (\chi\sin k\psi-\psi\sin k\chi)(\cos k\psi-\cos k\chi)>0,
\end{eqnarray*}
since both expressions in parentheses are positive. (The first one is
positive because $\sin(x)/x$ is decreasing on $(0,\pi/2)$.)

Therefore, $f(k)$ is increasing and we have
\begin{equation}\label{leq2}
\frac{1-\cos k\psi\cos k\chi}{\sin k\psi\sin k\chi}<
\frac{1-\cos \psi\cos \chi}{\sin \psi\sin \chi},
\quad 0<k<1.
\end{equation}
Since $l_1$ and $l_2$ are disjoint, $PQ$ is greater than in case of
parallel lines:
\begin{equation}\label{leq3}
\cosh PQ>\frac{1-\cos \psi\cos \chi}{\sin \psi\sin \chi}.
\end{equation}
Combining (\ref{leq1}), (\ref{leq2}), and (\ref{leq3}), we have
$\cosh PQ<\cosh PQ$. Contradiction.
\qed

\begin{lemma}\label{lemma2}
Let a triangle $\Delta=\Delta(\phi,\psi,\theta)$
in ${\Bbb H}^2$ have angles of
$\phi\leq \phi_0$, $\psi\leq\psi_0$, $\theta\leq\theta_0$
at vertices $N$, $L$, and $M$, respectively. Fix a real number
$c$ $(0<c\leq 1)$. Let $l=l(\phi,\psi,\theta)$ be the line
through $NM$ and let $d=d(\phi,\psi,\theta)$ be a half-line
from $L$ making an angle of $c\psi$ with $LM$ so that $d$ and
$\Delta$ lie in different half-planes with respect to $LM$.
Suppose $d$ and $l$ are disjoint for some $\phi=\phi_0$,
$\psi=\psi_0$, $\theta=\theta_0$. Then $d$ and $l$ are also
disjoint for any $\phi\leq\phi_0$, $\psi\leq\psi_0$,
$\theta\leq\theta_0$.
\end{lemma}

\noindent
{\it Proof.}
Fix a point $(\phi_1, \psi_1,\theta_1)$
in 3-dimensional parameter space such that
$\phi_1\leq\phi_0$, $\psi_1\leq\psi_0$,
$\theta_1\leq\theta_0$. Suppose that the lemma holds
at that point, i.e.,
$d=d(\phi_1, \psi_1,\theta_1)$
and
$l=l(\phi_1, \psi_1,\theta_1)$
are disjoint. It is sufficient to show that $d$ and $l$
are also disjoint at
$(\phi, \psi_1,\theta_1)$,
$(\phi_1, \psi,\theta_1)$, and
$(\phi_1, \psi_1, \theta)$, where
$\phi\leq\phi_1$, $\psi\leq\psi_1$,
$\theta\leq\theta_1$.

The first and the third cases are easy and true for any
$c>0$ that satisfies $c\psi_0<\pi$.

In case $(\phi_1,\psi,\theta_1)$, one can apply
Lemma~\ref{lemma_lines} with $l_2=NL$, $l_1=N'L'$,
$l_2'=d$, $l_1'=d'$, $k=c<1$ (see Figure~\ref{l1l2}b).
For $c=1$, the lines $d$ and $d'$
are disjoint, because they are images of $NL$ and $N'L'$
under reflection in $ML'$.
Hence, $d'$ and $l$ are also disjoint, since they lie on
different sides of $d$.
\qed

\subsection*{Polyhedra and links}

A plane divides ${\Bbb H}^3$ into two components;
we will call the closure of either of them a {\it half-space}
in ${\Bbb H}^3$.

A connected subset $P$ of ${\Bbb H}^3$ with non-empty interior
is said to be
a {\it (convex) polyhedron} if it is the intersection of a family
$\cal H$
of half-spaces with the property that each point of $P$
has a neighborhood meeting at most a finite number of
boundaries of elements of $\cal H$.

\medskip
\noindent
{\bf Definition.} Let $P$ be a polyhedron in ${\Bbb H}^3$ and let
$\partial P$ be its boundary in ${\Bbb H}^3$.
In (1)--(3) below we define the link for different
``boundary'' points of $P$ (cf.~\cite{EP}).

(1) Let $p\in\partial P$.
Let $S$ be a sphere in ${\Bbb H}^3$
with center $p$, whose radius is chosen small enough
so that it only meets faces of $P$ which contain $p$.
Such a sphere exists by the local finiteness property we claim in the
definition of a polyhedron.
There is a natural way to endow $S$ with a spherical
geometry identifying $S$ with ${\Bbb S}^2$ as follows.
Map conformally ${\Bbb H}^3$ onto the unit ball
${\Bbb B}^3=\lbrace x\in{\Bbb R}^3\ |\ |x|<1\rbrace$
so that $p$ goes to 0 and after that change the scale of
the sphere to be of radius~1.
The {\it link} of $p$ in $P$ is defined to be the image
of $S\cap P$ under the above identification
(it is well-defined up to isometry).

(2) Let $\overline{\partial P}$ be the closure of $\partial P$
in $\overline{\Bbb H}^3={\Bbb H}^3\cup\overline{\Bbb C}$.
Suppose
$\overline{\partial P}\backslash\partial P\not=\emptyset$,
and let $p\in\overline{\partial P}\backslash\partial P$.
Then $p\in\overline{\Bbb C}$ (i.e., it is an ideal point).
Let $S$ be a horosphere
centered at $p$ that only meets those faces of $P$
whose closures in $\overline{\Bbb H}^3$ contain $p$.
We can identify $S$ with Euclidean plane ${\Bbb E}^2$
using an isometry of ${\Bbb H}^3$ that sends $p$
to $\infty$.
The image of $S\cap P$ under such identification is called
the {\it link} of the ideal boundary point $p$ in $P$.
Note that such a link is defined up to similarity.

(3) Suppose that there exists a hyperbolic plane $S$ orthogonal to
some faces $F_1,\dots,F_t$ of $P$, and suppose that the other
faces of $P$ lie in the same open half-space which is
bounded by~$S$.
If $t\geq 3$ then we say that $S$ corresponds
to an {\it imaginary vertex} $p$ of $P$; and we define the {\it link}
of $p$ in $P$ to be $S\cap P$.

\medskip
Notice that the link of a proper (lying in ${\Bbb H}^3$),
ideal, or imaginary vertex $p$ in $P$ is
a spherical, Euclidean, or hyperbolic polygon (possibly
of infinite area), respectively.

\begin{figure}[htbp]
\centering
\includegraphics[width=8 cm]{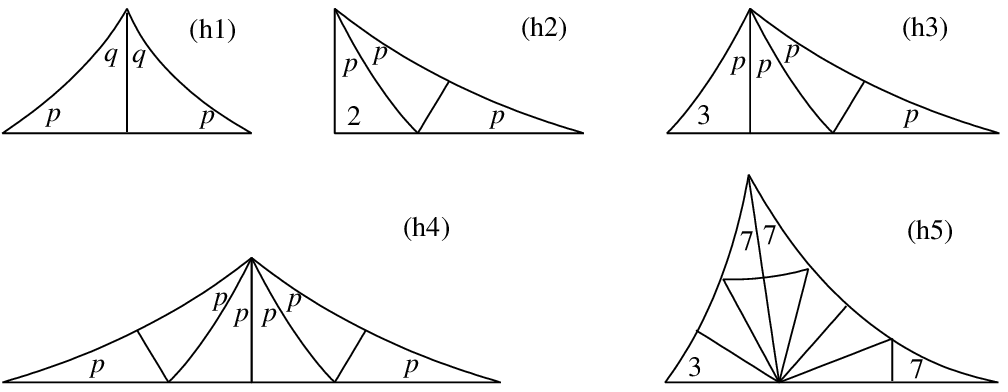}
\includegraphics[width=11 cm]{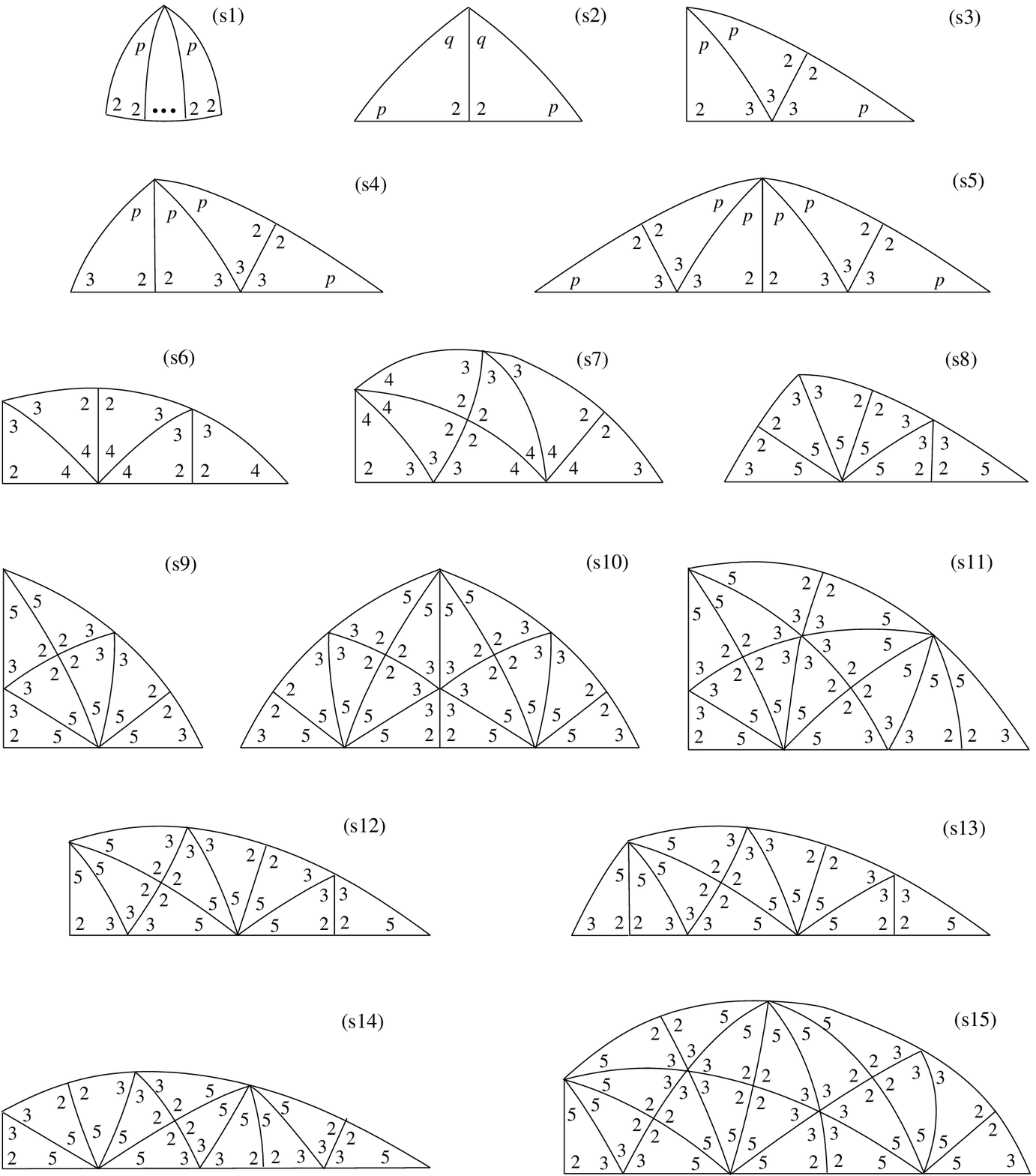}
\caption{All hyperbolic and spherical triangles with
two primitive and one non-primitive angles that generate a
discrete group}\label{triangles}
\end{figure}

The surface $S$ in the definition is called a {\it link surface} and
is orthogonal to all
faces of $P$ that meet $S$; hence the group generated by
reflections in these faces keeps $S$ invariant and can be
considered as the group of reflections in sides of a
spherical, Euclidean, or hyperbolic polygon
depending on the type of the vertex.

If such a polygon is a triangle with two primitive and one
non-primitive angle, we use Figure~\ref{triangles} that gives
a list of all such triangles in hyperbolic and spherical cases
which `generate' a discrete reflection group. The list of the hyperbolic
triangles first appeared in Knapp's paper~\cite{Kn}, then
Matelski gave a nice geometric proof~\cite{Ma}.
The list for the spherical case can be found in~\cite{F2}.
As for triangles with more than one
non-primitive angle reflections in sides of which generate a discrete group,
see~\cite[Lemma~2.1]{KS}
for a hyperbolic list and~\cite{F2} for a spherical one.

We denote a triangle with angles $\pi/p$, $\pi/q$, and $\pi/r$
in any of spaces ${\Bbb S}^2$, ${\Bbb E}^2$ or ${\Bbb H}^2$
by $(p,q,r)$. Note that the only Euclidean triangle with at least
one fractional angle that gives a discrete group is
$(6,6,3/2)$.

\medskip
\noindent
{\bf Definition.}
We define a {\it tetrahedron} $T$ to be a polyhedron which in the projective
ball model is the intersection of the hyperbolic space ${\Bbb H}^3$
with a usual Euclidean tetrahedron $T_E$ (possibly with vertices on
the sphere $\partial {\Bbb H}^3$ at infinity or beyond it) so that the
intersection of each edge of $T_E$ with ${\Bbb H}^3$ is non-empty.

Such a tetrahedron can be non-compact or have infinite volume, but the
links of its vertices are compact spherical, Euclidean, or hyperbolic
triangles.

Let a tetrahedron $T$ (possibly of infinite volume) have
dihedral angles $\pi/\lambda_1$, $\pi/\lambda_2$, $\pi/\lambda_3$
at some face and let $\pi/\mu_1$, $\pi/\mu_2$, $\pi/\mu_3$ be
dihedral angles of $T$ that are opposite to
$\pi/\lambda_1$, $\pi/\lambda_2$, $\pi/\lambda_3$, respectively.
We denote such a tetrahedron by
$T=[\lambda_1,\lambda_2,\lambda_3;\mu_1,\mu_2,\mu_3]$ and
the group generated by reflections in its faces by
$G_T$.

\section{Basic construction}

In this section we describe the main steps of the proof of
Theorem~A and prepare for a detailed description of the
list of all discrete groups.

The method of the proof is geometric and is based on using the
Poincar\'e theorem~\cite{EP}. However, the Poincar\'e theorem gives
conditions on the dihedral angles of a fundamental polyhedron for a group.
We would like to obtain conditions that depend only on
the generators $f$ and $g$ of $\Gamma$.
That is why we introduce auxiliary elements $h_i$, $i=1,\dots,4$,
in Theorem~A.

Each of $h_i$ is defined by two relations.
For example, for $h_1$ we have $h_1^2=gfg^{-1}f$ and
$(h_1f^{-1})^2=1$.
The first relation means that $h_1$ is a square root of $gfg^{-1}f$,
but there are in general two square roots of $gfg^{-1}f$ in $\PSL$.
In order to rewrite conditions used in the Poincar\'e theorem,
we need namely that one which is defined by
the second relation  $(h_1f^{-1})^2=1$.
The same is true for other elements $h_i$.

Existence of $h_1$ and $h_2$ is proved in Subsection~3
of this section
and existence of $h_3$ is proved in Section~4.1.
This proves the part (1)
of Theorem~A.
Existence of $h_4$ that appears in item~(2)$(vii)$ of Theorem~A
is proved in Section~4.2.

\bigskip
\noindent
{\it 3.1. Construction of $\Gamma^*$.}
We start with construction of a group $\Gamma^*$ containing $\Gamma$
as a subgroup of finite index. Such a group is discrete if and
only if $\Gamma$ is.

Let $f$ be a primitive elliptic element of odd order $n\geq 3$, $g$
be a hyperbolic element, and let their axes intersect
non-orthogonally. We denote elements and their axes by the same
letters when it does not lead to any confusion. Let $\omega$ be a
plane containing $f$ and $g$, and let $e$ be a half-turn with the
axis which is orthogonal to $\omega$ and passes through the point of
intersection of $f$ and $g$.

For the group $\Gamma=\langle f,g \rangle$ we define two finite
index extensions of it as follows:
$\widetilde{\Gamma}=\langle f,g,e\rangle$ and
$\Gamma^*=\langle f,g,e,R_\omega\rangle$
(we denote the reflection in a plane $\kappa$ by $R_\kappa$).
The groups $\Gamma$, $\widetilde\Gamma$,
and $\Gamma^*$ are either all discrete or all non-discrete.

\begin{figure}[htbp]
\centering
\begin{tabular}{cc}
\includegraphics[width=5 cm]{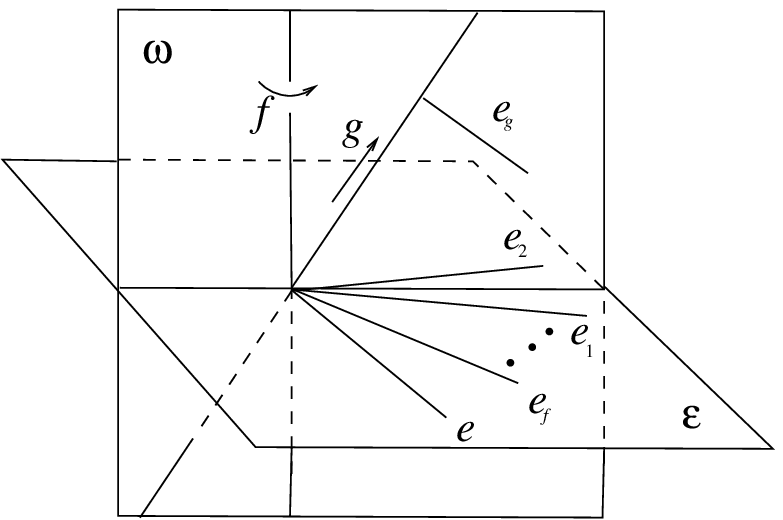}\qquad &
\qquad \includegraphics[width=4.2 cm]{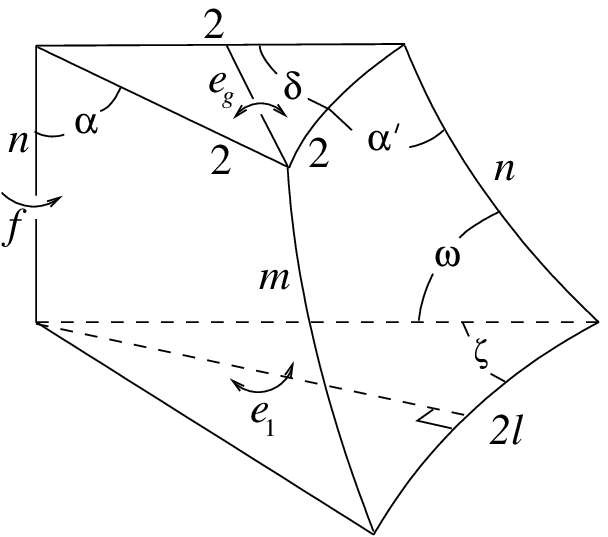}\\
(a)\qquad & \qquad(b)\\
\end{tabular}
\caption{}\label{planes}
\end{figure}

\bigskip
\noindent
{\it 3.2. New generators for $\Gamma^*$.}
Let $e_f$ and $e_g$ be half-turns such that $f=e_fe$ and $g=e_ge$.
Axes $e_f$ and $e$ lie in some plane, denote it by $\varepsilon$, and
intersect at an angle of $\pi/n$; $\varepsilon$ and $\omega$ are
mutually orthogonal; $e_g$ is orthogonal to $\omega$ and intersects
$g$, moreover, the distance between $e_g$ and $e$ is equal to half of
the translation length of $g$.

Consider $\varepsilon$ and $\langle e,e_f\rangle$
(see Figure~\ref{planes}a).
The group
contains elements $e$, $e_f=fe$, $f^2e$, \dots Each element $f^ke$,
$k=0,1,2,\dots$ is a half-turn with an axis lying in $\varepsilon$.
Denote $e_1=f^{(n-1)/2}e$ and $e_2=f^{(n+1)/2}e$. Note that
$\omega$ is a bisector of the two lines $e_1$ and $e_2$.

Let $\alpha$ be a hyperbolic plane such that $f=R_\omega R_\alpha$.
Then
$\Gamma^*=\langle e_1,e_g,R_\alpha, R_\omega\rangle$.

\medskip
\noindent
{\it 3.3. Construction of a fundamental polyhedron for $\Gamma^*$
corresponding to Theorem~A(2)({\it i}) and
existence of $h_1$ and $h_2$.}
First, note that
there exists a plane $\delta$ which is orthogonal to the planes
$\alpha$,
$\omega$, and $\alpha'=e_g(\alpha)$. The plane $\delta$ passes
through the
common perpendicular to $f$ and $e_g(f)$ orthogonally to $\omega$. It
is clear that $e_g\subset\delta$.

We need one more plane, denote it by $\zeta$, for the
construction of a fundamental polyhedron for $\Gamma^*$.
To construct $\zeta$ we use an auxiliary plane $\kappa$ that passes
through $e_1$ orthogonally to $\alpha'$. The plane $\zeta$ then
passes through $e_1$ orthogonally to $\kappa$.
(Note that in general $\zeta\neq\varepsilon$.)
In fact, the plane $\alpha'$ and the line $e_1$ can either
intersect, or be parallel, or disjoint. Note that if
$\zeta\cap \alpha'\not=\emptyset$ then
$e_1\perp(\zeta\cap \alpha')$.

Let $\P$ be the convex polyhedron bounded by $\alpha$, $\omega$,
$\alpha'$, $\delta$, and $\zeta$. Note that $\P$ can be compact
or non-compact (see Figure~\ref{planes}b, where $\P$ is drawn as
compact).

Consider the dihedral angles of $\P$. The angles between $\delta$ and
$\omega$, $\delta$ and $\alpha$, $\delta$ and $\alpha'$ are of
$\pi/2$; the angles formed by $\omega$ with $\alpha$ and $\alpha'$
are equal to $\pi/n$; the sum of the angles formed by $\zeta$ with
$\omega$ and $\alpha$ is equal to $\pi$. The planes $\alpha$ and
$\alpha'$ can either intersect or be parallel or disjoint. The
same is true for $\zeta$ and $\alpha'$. Denote the angle between
$\alpha$ and $\alpha'$ by $\pi/m$ and the angle between
$\alpha'$ and $\zeta$ by $\pi/(2l)$;
we use the symbols $\infty$ and $\overline\infty$ for $l$ or $m$ if
the corresponding planes are parallel or disjoint, respectively.
It is clear that if
\begin{equation}\label{condition-i}
\text{$m$ $(2/n+1/m<1)$ and $l$
$(l\geq 2)$ are integers, $\infty$, or
$\overline\infty$,}
\end{equation}
then $\P$ and elements $e_1$, $e_g$,
$R_\alpha$, $R_\omega$,
$R_{\alpha}'=e_gR_\alpha e_g$ satisfy the hypotheses of
the Poincar\'e
theorem, $\Gamma^*$ is discrete,
and $\P$ is its fundamental polyhedron.

Now we rewrite the condition~(\ref{condition-i}) via conditions on elements of
$\PSL$.
In fact, we wish to show that this sufficient condition for
discreteness of $\Gamma^*$ is equivalent to Item~({\it i})
of the theorem.

A natural choice of elements $R'_\alpha R_\alpha$ and
$R''_\alpha R'_\alpha$, where $R''_\alpha=e_1 R_\alpha e_1$, does
not provide enough information about $\P$. (Even if $R'_\alpha
R_\alpha$ and $R''_\alpha R'_\alpha$ are primitive elliptic
elements, the corresponding dihedral angles of $\P$ can be
obtuse. We refer the reader to \cite[pp.~257--258]{KK} for more
detailed explanation.)
Therefore, we choose the following square roots instead.

The first element we define is $h_1=R_\xi R_\alpha=R'_\alpha R_\xi$,
where $\xi$ is the bisector
of $\alpha$ and $\alpha'$ passing through $e_g$. Then
$$
h_1^2=R'_\alpha R_\alpha=R'_\alpha R_\omega R_\omega R_\alpha=f'f,
$$
where
\begin{equation}\label{fprime}
f'=R'_\alpha R_\omega=e_gR_\alpha e_gR_\omega
=e_gR_\alpha R_\omega e_g=e_gf^{-1}e_g=e_gefee_g=gfg^{-1}.
\end{equation}
Since $\xi$ is orthogonal to $\omega$, $(R_\xi R_\omega)^2=1$. On the
other hand,
$$
R_\xi R_\omega=h_1R_\alpha R_\omega=h_1f^{-1}.
$$
The above equations imply two conditions
$$
h_1^2=gfg^{-1}f \qquad {\rm and} \qquad (h_1f^{-1})^2=1
$$
that uniquely determine $h_1$ as an element of $\PSL$. Moreover,
$h_1$ is a primitive elliptic element of even order $2m$
$(2/n+1/m<1)$ if and only if the dihedral angle of $\P$ at the
edge $\alpha\cap\alpha'$ is equal to $\pi/m$, $m\in {\Bbb Z}$
\ $(2/n+1/m<1)$;
$h_1$ is parabolic (hyperbolic) if and only if $\alpha$ and
$\alpha'$ are parallel (disjoint, respectively).

Similarly, we define the second element as $h_2=R_\zeta
R_\alpha'=R_\alpha''R_\zeta$. Then
\begin{equation}\label{h22}
h_2^2=R_\alpha''R_\alpha'=R_\alpha''R_\omega R_\omega R_\alpha'=
t(f')^{-1}=tgf^{-1}g^{-1},
\end{equation}
where
$$
t=R_\alpha''R_\omega=e_1R_\alpha'e_1R_\omega=
e_1R_\alpha'R_\alpha e_1=e_1h_1^2e_1.
$$
Furthermore, since $e_1=f^{(n-1)/2}e$ and $h_1^2=gfg^{-1}f$,
\begin{equation}\label{t}
t=f^{(n-1)/2}egfg^{-1}f^{(n+1)/2}e=
f^{(n-1)/2}g^{-1}f^{-1}gf^{-(n+1)/2}.
\end{equation}
From (\ref{h22}) and (\ref{t}) we have
$h_2^2=f^{(n-1)/2}g^{-1}f^{-1}gf^{-(n+1)/2}gf^{-1}g^{-1}$.

We now need to choose the square root of
$R_\alpha''R_\alpha'$ that corresponds to the dihedral angle of $\P$
made by $\alpha'$ and $\alpha''$. The required $h_2$ is such that
$h_2f'=h_2gfg^{-1}$ is an elliptic element whose axis intersects $f$.

With this choice, $h_2$ is a primitive elliptic element of
even order $2l$, $l\geq 2$, if and only if the dihedral angle between
$\alpha'$ and $\alpha''$ is $\pi/l$, $l\in {\Bbb Z}$ $(l\geq 2)$;
$h_2$ is
parabolic (hyperbolic) if and only if $\alpha'$ and $\alpha''$ are
parallel (disjoint, respectively).

Therefore, the condition~(\ref{condition-i}) is equivalent to
{\setlength\arraycolsep{0pt}
\begin{eqnarray*}
&&\text{$h_1$ is hyperbolic, parabolic or a primitive
elliptic element of even}\nonumber\\
&&\text{order $2m$ \ $(2/n+1/m<1)$ \ and $h_2$ is  hyperbolic,
parabolic or}\\
&&\text{a~primitive elliptic
element of even order $2l$ $(l\geq 2)$,}\nonumber
\end{eqnarray*}
}which is item~(2)$(i)$ of  Theorem~A.

\medskip
\noindent
{\it 3.4.
Assume that the condition~(\ref{condition-i}) does not hold.}
This means that $\P$ itself is not a fundamental polyhedron for
$\Gamma^*$. Then  discrete groups may appear in one of the following cases:
\begin{itemize}
\item[(1)] $m$ is fractional, $1/m+2/n<1$;
\item[(2)] $m$ is an integer, $1/m+2/n<1$, and $l$ is fractional,
$l>1$;
\item[(3)] $m\in\{\infty,\overline\infty\}$ and $l$ is fractional, $l>1$.
\end{itemize}

These cases require a further consideration and
in Sections~4--6 we completely determine the list of discrete groups
in each case.

Namely, in each of the cases (1)--(3) we assume first
that $\Gamma^*$ is discrete
to exclude {\it some} groups that cannot be discrete,
that is, we get necessary conditions for discreteness.
It turns out that {\it all} remaining groups are discrete,
that is the necessary conditions are also sufficient.
Such discrete groups are listed in items~(2)$(ii)$--$(vii)$ of
Theorem~A. In order to prove sufficiency, we forget that
$\Gamma^*$ was assumed to be discrete and
for each such group we give a polyhedron
and generators of $\Gamma^*$ satisfying the hypotheses of the
Poincar\'e theorem. Hence, $\Gamma^*$, and therefore $\Gamma$, is discrete
and the sufficient part of (2) in Theorem~A is also proved.

\section {$m$ is fractional, $1/m+2/n<1$}

We have assumed that
$\Gamma^*=\langle e_1,e_g,R_\alpha,R_\omega\rangle$
is discrete, so each of its subgroups is also discrete. Hence
$\langle R_\omega,R_\alpha,R'_\alpha\rangle$ is discrete. Notice that
$\langle R_\omega,R_\alpha,R'_\alpha\rangle$ acts as a group of
reflections in the sides of a hyperbolic triangle $(n,n,m)$ which
is
the upper face of $\P$ (see Figure~\ref{planes}b). From the list of
all triangles with two
primitive angles that give a discrete group (see
Figure~\ref{triangles}),
we have that for $m$ fractional either
\begin{itemize}
\item[(a)] $\widetilde m=2m$ is odd; or
\item[(b)] $4m=n$.
\end{itemize}
In both cases $\Gamma^*$ contains the reflection $R_\xi$, where $\xi$
is the  plane that bisects the dihedral
angle of $\P$ at the edge $\alpha\cap \alpha'$.
Moreover, $R_\delta$ also belongs to $\Gamma^*$, because $\xi$ passes
through $e_g$ and $e_g=R_\xi R_\delta$.

\subsubsection*{4.1. $\widetilde m=2m$ is odd}
That is
\begin{equation}\label{condition-v-1}
\text{$h_1$ is a primitive elliptic
element of odd
order $\widetilde m$, $1/n+1/\widetilde m<1/2$.}
\end{equation}

\begin{figure}[htbp]
\centering
\includegraphics[width=4 cm]{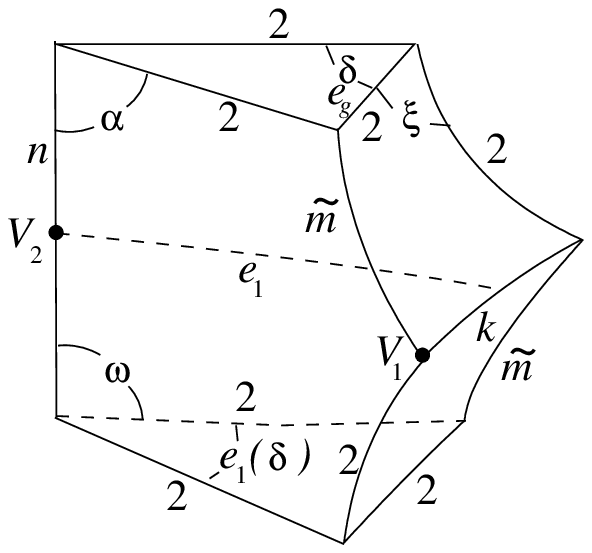}
\caption{Polyhedron $\P_1$}\label{polyhedron_p1}
\end{figure}

\noindent
Let $\P_1$ be a polyhedron bounded by $\omega$, $\alpha$,
$\delta$, $e_1(\delta)$, $\xi$, and $e_1(\xi)$. Let $\pi/k$ be the
dihedral angle of $\P_1$ at the edge $\xi\cap e_1(\xi)$,
$k\in(1,\infty)\cup\{\infty,\overline\infty\}$. The other angles of
$\P_1$ are submultiples of $\pi$ (see
Figure~\ref{polyhedron_p1}).

\noindent
If
\begin{equation}\label{condition-v-2}
\text{$k>1$ is an integer, $\infty$, or $\overline\infty$}
\end{equation}
then
$\Gamma^*$ is discrete. Its fundamental polyhedron is a half of
$\P_1$.

Again, rewrite the condition~(\ref{condition-v-2}) as
a condition  on elements of $\PSL$.
Let $h_3=R_\beta R_\xi$, where $\beta$ is the bisector of $\xi$ and
$e_1(\xi)$ passing through $e_1$.
Then
\begin{eqnarray*}\label{h32}
h_3^2&=&(R_\beta R_\xi)^2=(e_1R_\xi e_1)R_\xi=
e_1R_\xi e_1R_\omega R_\omega R_\xi=
e_1R_\xi R_\alpha e_1h_1f^{-1}\\
&=&e_1h_1e_1fh_1^{-1}
=f^{(n-1)/2}eh_1ef^{-(n-1)/2}fh_1^{-1}.
\end{eqnarray*}
Note that since $e_gh_1e_g=h_1^{-1}$, we have $eh_1e=g^{-1}h_1^{-1}g$
and,
therefore,
$$
h_3^2=f^{(n-1)/2}g^{-1}h_1^{-1}gf^{-(n-3)/2}h_1^{-1}.
$$
We have $h_3=R_\beta R_\xi$ if and only if $h_3h_1$ is an elliptic element
whose axis intersects the axis of $f$.

The element $h_3$ is a primitive elliptic of even order $2k$,
$k\geq 2$, if and only if the dihedral angle of $\P_1$
between $\xi$ and $e_1(\xi)$ is $\pi/k$, $k\in {\Bbb Z}$,
$k\geq 2$; $h_3$ is
parabolic (hyperbolic) if and only if $\xi$ and $e_1(\xi)$ are
parallel (disjoint, respectively).

So we have that
the conditions~(\ref{condition-v-1}) and~(\ref{condition-v-2})
are equivalent to Item~$(v)$ of Theorem~A.

Assume that $k$ is fractional. Then there are additional reflections
in
$\Gamma^*$ in planes passing through the edge $\xi\cap e_1(\xi)$
which
decompose $\P_1$. Consider the link of the vertex
$V_1=\alpha\cap\xi\cap e_1(\xi)$ (see Figure~\ref{polyhedron_p1}). It is either
a spherical, Euclidean,
or hyperbolic triangle $(2,\widetilde m,k)$, $\widetilde m$ is odd,
depending on whether $V_1$ is a
proper, ideal, or imaginary vertex, respectively.
Since the corresponding link surface $S_1$
is invariant under
$\langle R_\alpha,R_\xi,R'_\xi\rangle$, where $R'_\xi=e_1R_\xi e_1$,
the group naturally
acts on $S_1$ as a group of reflections in the sides of the triangle
$(2,\widetilde m,k)$.
There are no discrete
groups with $\widetilde m$ odd and $k$ fractional if
$(2,\widetilde m,k)$ is
Euclidean. As for spherical and hyperbolic triangles, there are the
following types of the link of $V_1$:
\begin{itemize}
\item spherical triangles (s3), (s7), (s9), (s11), (s12), (s14),
(s15);
\item hyperbolic triangle (h2) in Figure~\ref{triangles}.
\end{itemize}

\begin{lemma}\label{lemma1}
If $\widetilde m$ is odd and the bisector $\beta$ of the dihedral
angle of $\P_1$ as above
formed by $\xi$ and $e_1(\xi)$ is the plane of a reflection belonging to
the discrete group $\Gamma^*$, then either $n=3$ and
$\widetilde m\geq 7$, or
$n=5$ and $\widetilde m\geq 5$.
\end{lemma}

\medskip
\noindent
{\it Proof.} Let $\beta$ intersect $\alpha$ at an angle of
$\pi/p$,
where
$2<p<\widetilde m$ ($\beta$ lies
between $\xi$ and $e_1(\xi)$,
$p$ is not necessary an integer).
Then $R_\beta R_\alpha$ is an elliptic element of order greater
than $2$. Moreover, the axis of this element meets $f$ at
$V_2=f\cap e_1$. We have two elliptic elements $f$ and
$R_\beta R_\alpha$ with different intersecting axes whose orders
are greater than 2. Since these elements belong to a discrete
group, the orders are at most 5. So, for odd $n$ we have $n=3$ or
$n=5$.
It follows from the
upper face of $\P_1$ that $1/n+1/\widetilde m<1/2$.
Thus, $n=3$ implies
$\widetilde m\geq 7$ and $n=5$ implies
$\widetilde m\geq 5$.\qed

\medskip

In {\bf Cases (s9)} and {\bf (s15)} $\widetilde m$ should be
3, but this contradicts Lemma~\ref{lemma1} and these cases disappear
(the group is not discrete).
We will use Lemma~\ref{lemma1} also to specify $n$ and
$\widetilde m$ in Cases~(s3) and (s14).

\medskip 

Let us describe the construction for {\bf Case~(s3)} in detail. All
the other decompositions are made in a similar manner and we will
give only a short description for each of them.

By Lemma~\ref{lemma1}, we have $n=5$ and $\widetilde m=5$.
Let us first consider a non-compact polyhedron $\widetilde\P_1$
bounded by $\omega$, $\alpha$, $\xi$ and $e_1(\xi)$. The reflection
planes $\delta$ and $e_1(\delta)$ will be added later.

\begin{figure}[htbp]
\centering
\begin{tabular}{cc}
\includegraphics[width=4 cm]{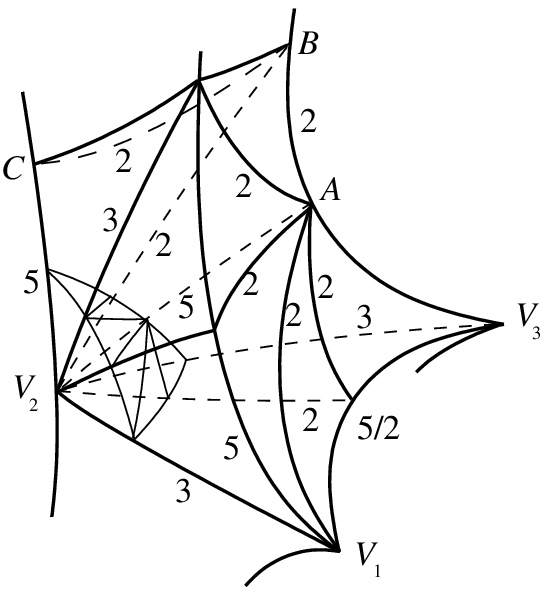}\quad  &
\quad \includegraphics[width=4 cm]{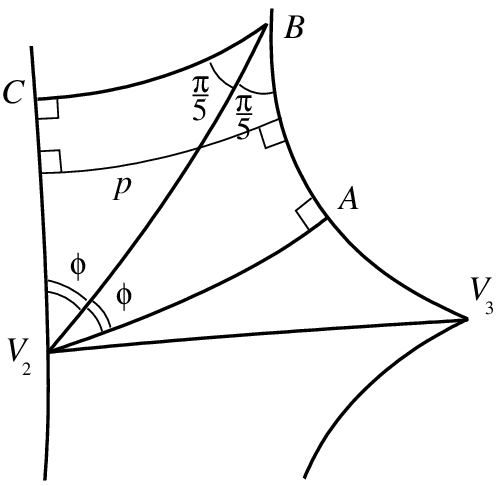}
\\
(a)\qquad & \qquad(b)\\
\end{tabular}
\caption{}\label{p1_s3}
\end{figure}

Since the bisector $\beta$ is a reflection plane, there exists a
reflection plane $\eta$ passing through $e_1$ orthogonally to $\beta$.
From the link of $V_1$, there exists a reflection plane $\eta_1$
through $\alpha\cap\beta$ that makes an angle of $\pi/3$ with $\alpha$
(see Figure~\ref{p1_s3}a).
The link of the vertex $V_2$
formed by $\alpha$, $\omega$ and $\beta$ is a spherical triangle
$(5,3,3/2)$, which can be decomposed
into smaller triangles only as (s8). Then
$\widetilde\P_1$, and therefore ${\Bbb H}^3$, is decomposed into tetrahedra
$T=[2,2,5;2,3,5]$ (see the last paragraph in Section~2, which
explains this notation).
Each of $T$ is a fundamental polyhedron for the
group $\langle R_\omega,R_\alpha,R_\xi,e_1\rangle=G_T < \Gamma^*$.

Now we determine the position of $\delta$. Consider the face of
$\widetilde\P_1$ lying in $\omega$. In Figure~\ref{p1_s3}b we have
drawn traces of the reflection planes in $\omega$. (Note that
not all of these traces are reflection lines in $\omega$.)
The angle at $B$ of the quadrilateral $ABCV_2$ is $2\pi/5$. From the
link of $V_2$,
$$
\cos\phi=\frac{\cos(\pi/3)+\cos(\pi/2)\cos(\pi/5)}
{\sin(\pi/2)\sin(\pi/5)}=
\frac{1}{2\sin(\pi/5)}>\frac{1}{2\sin(\pi/4)}=\frac 1{\sqrt{2}}.
$$
That is $2\phi<\pi/2$. Thus, the common perpendicular $p$ to $CV_2=f$
and $AB=\omega\cap\xi$ lies inside $ABCV_2$. Since $\delta$ passes
through $p$, $\delta$ intersects the interior of one of the tetrahedra
$T$. By \cite{DM}, $\Gamma^*=\langle G_T,R_\delta\rangle$
is not discrete.

\medskip
In {\bf Case~(s14)}, $n=5$ and $\widetilde m=5$, the space is
decomposed into tetrahedra $T=[2,2,3;2,5,3]$. By the same reason as in
Case~(s3), $\Gamma^*$ is not discrete.
So, it remains to
consider only Cases~(s7), (s11), (s12) without the central plane and
the hyperbolic case~(h2).

\begin{figure}[htbp]
\centering
\begin{tabular}{cc}
\includegraphics[width=3.5 cm]{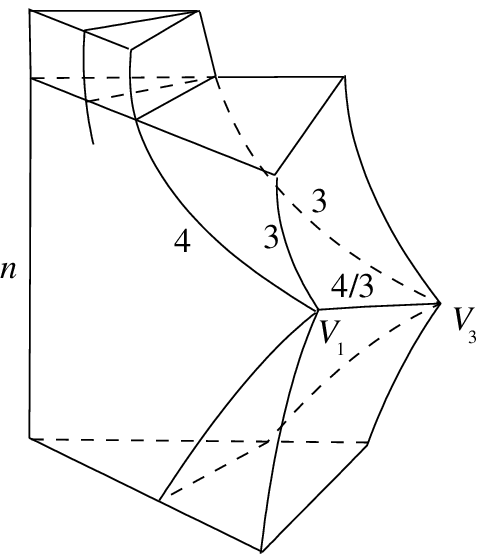}\qquad &
\qquad \includegraphics[width=4 cm]{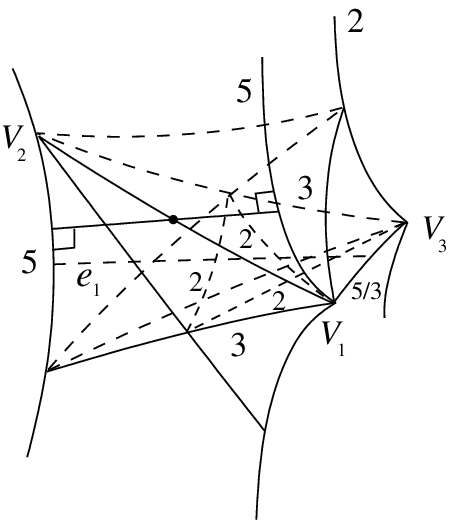}\\
(a)\qquad & \qquad(b)\\
\end{tabular}
\caption{}\label{p1_s7_s12}
\end{figure}

\medskip
{\bf Cases (s7)} and {\bf (s11)}.
Here $k=p/3$, $p=4$ or $5$, $\widetilde m=3$,
$n\geq 7$.
There are two planes through $\xi\cap e_1(\xi)$. Consider one of them
(denote it by $\eta$). From the links of $V_1$ and $V_3$ (see
Figure~\ref{p1_s7_s12}a for Case~(s7)) one can easily find the
dihedral angles which
$\eta$ forms with $\alpha$ and $\omega$. Since $n\geq 7$ and
$\Gamma^*$ is supposed to be discrete, the vertex formed by $\alpha$,
$\omega$, and $\eta$ is imaginary and its link is a hyperbolic
triangle $\Delta=(n,4,3/2)$ for Case~(s7) and
$\Delta=(n,3,5/3)$ for (s11);
$\langle R_\alpha,R_\eta, R_\omega\rangle$ acts
as a group of reflections in the sides of $\Delta$. But
there
are no such discrete groups. Thus, $\Gamma^*$ is not discrete.

\medskip
{\bf Case (s12).} Here $k=5/3$, $\widetilde m=5$, $n\geq 5$.
Consider the polyhedron of infinite volume bounded by $\alpha$,
$\omega$, $\xi$, and $e_1(\xi)$.
The link of the vertex $V_2$ is a spherical, Euclidean, or hyperbolic
triangle $(n,5,3/2)$. Since the group generated by the reflections in
its sides is discrete, $n=5$.
It is easy to find the decomposition of the space
into
tetrahedra $T=[2,3,5;2,3,2]$ (see Figure~\ref{p1_s7_s12}b).
One can see that the common perpendicular
$p$ to the lines $f$ and $\alpha\cap\xi$ intersects $V_1V_2$ in the
middle. Thus, the plane $\delta$ (which is orthogonal to $\alpha$,
$\omega$, and $\xi$) passes through $p$ and cuts $T$ into two parts.
But $\langle G_T,R_\delta\rangle<\Gamma^*$ is not discrete by
\cite{DM}.

\medskip
{\bf Case (h2).} Here
$k=\widetilde m/2$, $\widetilde m\geq 7$, $(\widetilde m,2)=1$, $n=3,5$.

Assume first that $n=3$. Consider the link of $V_1$. The
bisector
$\beta$ of the dihedral angle formed by $\xi$ and $e_1(\xi)$ passes
through $e_1$. Since $R_\beta\in\Gamma^*$, there is a plane $\eta$
passing through $e_1$ orthogonally to $\beta$ and
$R_\eta\in\Gamma^*$. Note that $\eta$ is orthogonal to $\xi$. The
link of the vertex $V_2$ formed by $\alpha$, $\omega$, and $\beta$ is
a spherical triangle $(3,3,3/2)$ and is divided uniquely into 4
triangles $(2,3,4)$ by reflection planes.

\begin{figure}[htbp]
\centering
\begin{tabular}{cc}
\includegraphics[width=4 cm]{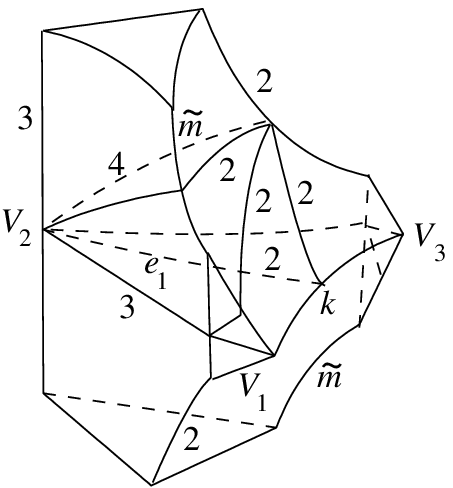}\qquad &
\qquad \includegraphics[width=4.2 cm]{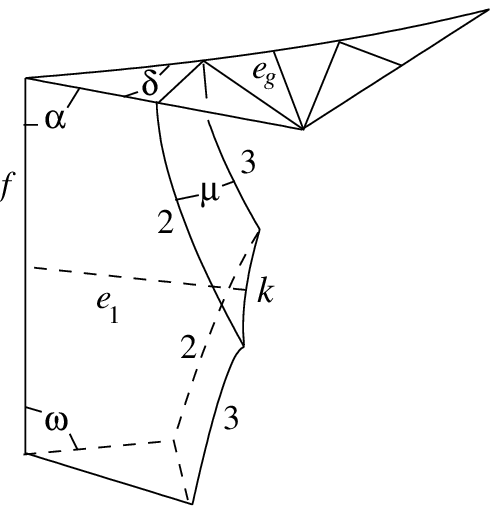}\\
(a)\qquad &\qquad (b)\\
\end{tabular}
\caption{Discrete groups for Items~({\it vi}) and~({\it vii})
of Theorem~A.}\label{p1_hyp}
\end{figure}

Let $\widetilde T$ be the prism bounded by $\omega$, $\beta$,
$\eta$, $\xi$, and the link surface $S$ of $V_3$ (see
Figure~\ref{p1_hyp}a). Since $R_\delta\in\Gamma^*$, it is clear that
the hyperbolic plane $S$
becomes a reflection plane,
$\Gamma^*=\langle
R_\omega,R_\beta,R_\eta,R_\xi,R_S\rangle$, and $\widetilde T$ is a fundamental
polyhedron for $\Gamma^*$. This is item~(2)({\it vi}) of the theorem.

Now assume that $n=5$. In this case the polyhedron $\widetilde\P_1$ is
decomposed into infinite volume tetrahedra $T=[2,2,5;2,3,\widetilde m]$.
Similarly to Case (s3), the plane $\delta$ intersects the interior of
$T$ so that $\Gamma^*$ is not discrete.

\subsubsection*{4.2. $m=n/4$, $n\geq 7$}
Consider the polyhedron $\P_2$ bounded by
$\alpha$, $\omega$, $\delta$, $e_1(\delta)$, $\mu$, and
$e_1(\mu)$ (see Figure~\ref{p1_hyp}b for determining the plane
$\mu$ which is orthogonal to $\delta$;
$\P_2$ is drawn as a compact). Denote the dihedral
angle of $\P_2$ formed by $\mu$ and $e_1(\mu)$ by
$\pi/k$. Note that
$\P_2$ is congruent to the polyhedron $\P_1$ with
$\widetilde m=3$, see Subsection~(a) above,
but we see that there are no discrete groups in
that case
with fractional $k$.
However, if
\begin{equation}
\text{$k>1$ is an integer, $\infty$ or
$\overline\infty$,}\label{condition-vii}
\end{equation}
then $\Gamma^*$ is discrete.

We show how to get Item~$(vii)$ of Theorem~A.
The procedure is analogous to the
usual one. First of all, note that now the element $h_1$ defined
in Section~3
is a rotation through $\pi/m=4\pi/n$
and we take a primitive elliptic $\tilde h_1$ of
order $n$ so that $h_1=\tilde h_1^2$. Further, define $h_4$ as
follows:
\begin{eqnarray*}
h_4^2&=&(e_1R_\mu e_1)R_\mu
=(e_1R_\mu e_1R_\omega)(R_\omega R_\mu)\\
&=&(e_1f\tilde h_1^{-1}fe_1)(\tilde h_1 f^{-1})=
f^{(n-3)/2}g^{-1}\tilde h_1gf^{-(n+1)/2}\tilde h_1f^{-1}.
\end{eqnarray*}
Note that to obtain the last equality we have used the obvious fact that
$e_g\tilde h_1^{-1}e_g=\tilde h_1$.

In order to define $h_4$ uniquely we need the last condition:
$h_4f\tilde h_1^{-1}$ must be an elliptic element whose axis
intersects $f$.
Now we see that the conditions ``$m=n/4$, $n\geq 7$''
and~(\ref{condition-vii}) are equivalent to
Item~({\it vii}) of the theorem.

\section{$m$ is an integer $(2/n+1/m<1)$ and $l$ is fractional
$(l>1)$}

Let $\widetilde\P$ be a non-compact polyhedron bounded by
$\alpha$, $\omega$, $\alpha'$, and
$\alpha''=e_1(\alpha')$. The intersection of $\alpha$, $\alpha'$, and
$\alpha''$ forms a vertex $V_1$. Its link is
either a spherical, Euclidean, or hyperbolic triangle $(n,m,l)$.
The subgroup $\langle R_\alpha,R_{\alpha}',R_{\alpha}'\rangle$ of
the discrete group $\Gamma^*$ is also discrete and
acts as the group
generated by reflections in the sides of the triangle $(n,m,l)$,
where $n$ is odd, $m$ is an integer, and $l$ is fractional.

There are no discrete groups if $V_1$ is ideal. As for proper or
imaginary $V_1$, using Figure~\ref{triangles} and
taking into account the fact that $2/n+1/m<1$,
we have the following list to consider (we also indicate
triangles
from Figure~\ref{triangles} that correspond to
$V_1$):
\begin{itemize}
\item[(1)] $n\geq 5$, $m=n$, $l=k/2$, $k$ is odd, $k\geq 3$ (h1, s2);
\item[(2)] $n\geq 5$, $m=2$, $l=n/2$ (h2, s3);
\item[(3)] $n=3$, $m\geq 4$, $(m,3)=1$, $l=m/3$ (h3, s4);
\item[(4)] $n\geq 5$, $(n,3)=1$, $m=3$, $l=n/3$ (h3, s4);
\item[(5)] $n\geq 5$, $m=n$, $l=n/4$ (h4, s5);
\item[(6)] $n=3$, $m=5$, $l=3/2$ (s8);
\item[(7)] $n=5$, $m=3$, $l=3/2$ (s8);
\item[(8)] $n=5$, $m=2$, $l=5/3$ (s12);
\item[(9)] $n=3$, $m=5$, $l=5/4$ (s13);
\item[(10)] $n=5$, $m=3$, $l=5/4$ (s13);
\item[(11)] $n=5$, $m=2$, $l=3/2$ (s14);
\item[(12)] $n=3$, $m=7$, $l=7/2$ (h5);
\item[(13)] $n=7$, $m=3$, $l=7/2$ (h5).
\end{itemize}

For each case we proceed as follows. We start with reflections in
$\alpha$, $\omega$, $\alpha'$, and $\alpha''$.
Using the link of $V_1$ we
find additional reflection planes and
dihedral angles that they make with
$\alpha$ and $\alpha'$.
Then using new links we
find new reflection planes and construct a
decomposition of $\widetilde\P$ into subpolyhedra (usually
tetrahedra) $T$ with primitive angles, or show that the group is not
discrete.
If the bisector $\zeta$ of the dihedral angle of
$\widetilde\P$ formed by $\alpha'$ and $\alpha''$ is a
reflection plane in $\Gamma^*$, then since it passes through
$e_1$, there is another reflection plane $\eta$ through $e_1$ that
is orthogonal to $\zeta$. We include $\eta$ in the decomposition.
If the decomposition gives a discrete group $G_T$
generated by reflections in the faces of $T$, we find the positions of
$e_1$ and $e_g$ and determine whether the group
$\langle G_T,e_1,e_g\rangle=\Gamma^*$ is discrete.
When construction of such a decomposition is not difficult we give
only a description of $T$.

\medskip
\noindent
{\bf Remark~1.}
If two axes of elliptic elements intersect non-orthogonally in a
discrete group, then the orders of both elements are at most $5$.
In particular, if $\Gamma^*$ is discrete and
the bisector $\zeta$ is a reflection plane in
$\Gamma^*$ so that $\zeta$ meets $\alpha$ non-orthogonally, then
$n=3$ or $5$ and the order of $R_\alpha R_\zeta$ is $3,4$, or $5$.

\begin{figure}[htbp]
\centering
\includegraphics[width=3 cm]{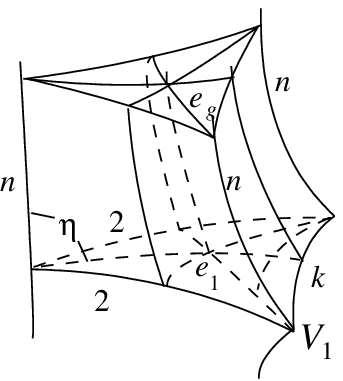}
\caption{}\label{ptf_1}
\end{figure}

\medskip
In {\bf Case (1)}, $n\geq 5$, $m=n$, $l=k/2$, $k$ is odd, $k\geq 3$.
The bisector
$\zeta$ is the only reflection plane passing through
$\alpha'\cap\alpha''$.
It is orthogonal to both $\alpha$ and $\omega$ and so to
$f=\alpha\cap\omega$.
Since $\zeta$ passes through $e_1$, there is
a reflection plane $\eta$ passing through $e_1$ orthogonally to
$\zeta$.
Since $\zeta$ is orthogonal to $f$, the plane $\eta$ contains $f$
and  is the bisector of the dihedral angle formed by
$\alpha$ and $\omega$.
Since also $m=n$, the polyhedron is symmetric with respect to
$\eta$. Besides $R_\eta$, $\Gamma^*$ contains
reflections $e_gR_\eta e_g$ and $R_\eta(e_gR_\eta e_g)R_\eta$ (see
Figure~\ref{ptf_1} for proper $V_1$). Note that the plane $\eta'$ of the last
reflection passes through $\alpha\cap\alpha'$ and makes an angle
of $\pi/(2n)$ with $\alpha$.

Now we show that $V_1$ cannot be proper, that is, $n\neq 5$.
Indeed, if $n=5$ then there are two rotations $R_\alpha R'_\eta$
of order $2n=10$ and $R'_\alpha R_\zeta$ of order $k$
($k$ is odd) with axes intersecting at $V_1$, which is impossible
in discrete groups.

$V_1$ cannot be imaginary as well. Indeed, from~\cite{KS} it follows
that the hyperbolic link of $V_1$ is a triangle $(n,n,n)$, but this
means that $k=2n$ is even.

So, Case (1) does not give discrete groups.

\begin{figure}[htbp]
\centering
\begin{tabular}{cc}
\includegraphics[width=3.5 cm]{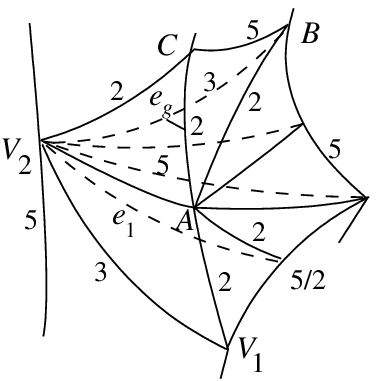} \qquad &
\qquad\includegraphics[width=3.9 cm]{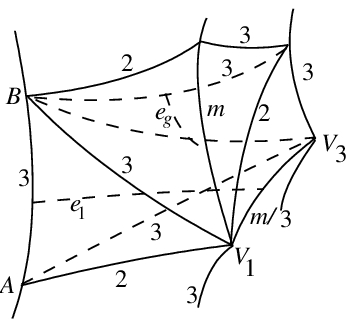}\\
(a)\qquad & \qquad (b)\\
\includegraphics[width=3.8 cm]{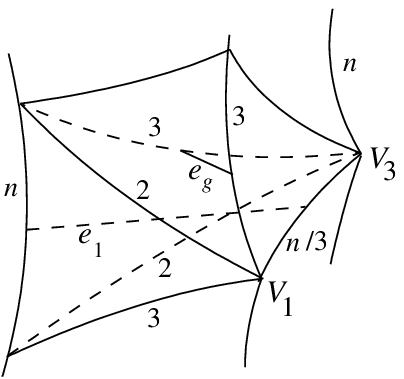} \qquad &
\qquad\includegraphics[width=4.8 cm]{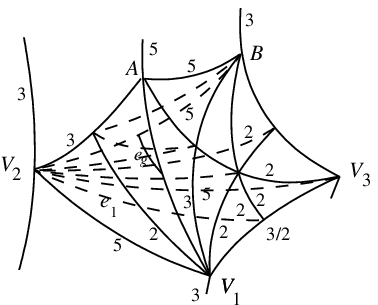}\\
(c)\qquad & \qquad (d)\\
\includegraphics[width=3.7 cm]{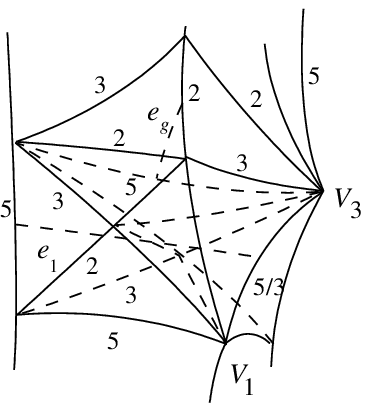} \qquad &
\qquad\includegraphics[width=3.2 cm]{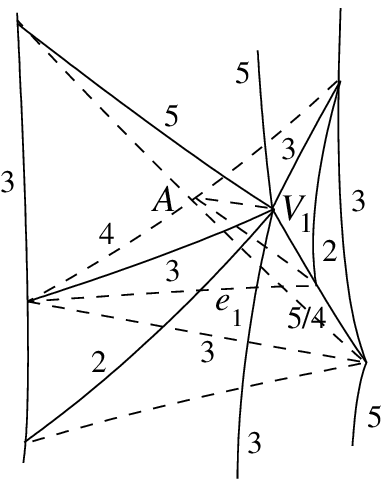}\\
(e)\qquad & \qquad (f)\\
\includegraphics[width=3.3 cm]{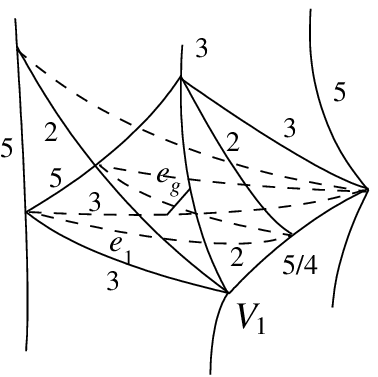} \qquad &
\qquad\includegraphics[width=3.0 cm]{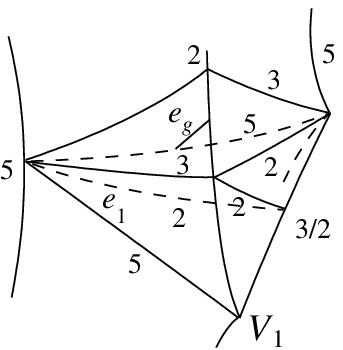}\\
(g)\qquad & \qquad (h)\\
\end{tabular}
\caption{Discrete groups for Item~({\it ii}) of Theorem~A.}\label{ptf_discr}
\end{figure}

\medskip
In {\bf Case (2)}, $n\geq 5$, $m=2$, $l=n/2$. The link of the vertex
$V_2$ formed by $\alpha$, $\omega$, and $\zeta$ is a
spherical triangle $(n,3,3/2)$. By Remark~1,
$n=5$. In this case $V_1$ is
proper and $\widetilde\P$ is decomposed into tetrahedra
$T=[2,3,5;2,2,5]$  so that $e_1$ is one of its edges
(see Figure~\ref{ptf_discr}a).

One can see that the line passing through the midpoints of the edges
$AC$ and $BV_2$ in the tetrahedron $ABCV_2$ is orthogonal to $\omega$
and to $V_1C$ and thus coincides with $e_g$. By \cite{DM}, the
group $\langle G_T, e_g\rangle=\Gamma^*$ is discrete and  a half
of $T$ is its fundamental polyhedron. This is a group from the
list in item~(2)({\it ii}) of the theorem.

\medskip
In {\bf Case (3)}, $n=3$, $m\geq 4$, $(m,3)=1$, $l=m/3$,
$\widetilde\P$ is decomposed into tetrahedra (possibly of
infinite volume) $T=[2,3,m;2,3,3]$ (see Figure~\ref{ptf_discr}b for
$m=4,5$). By construction, $e_1$ is orthogonal to $f$ and to
$V_1V_3$ and passes through the midpoint of $AB$. It is clear that
$\langle G_T, e_1\rangle$ is discrete. Note that $e_g$
is conjugate to $e_1$ and $\langle G_T, e_1\rangle=\Gamma^*$.
This is one of the groups from item~(2)({\it ii}) of the theorem.

\medskip
In {\bf Case (4)}, $n\geq 5$, $(n,3)=1$, $m=3$, $l=n/3$,
$\widetilde\P$ is decomposed into tetrahedra (possibly of
infinite volume) $T=[2,3,n;2,3,n]$ (see Figure~\ref{ptf_discr}c for
$n=5$). One can see that $e_1$ is orthogonal to the opposite edges
of $T$ with the dihedral angles of $\pi/n$, and $e_g$ is
orthogonal to the opposite edges with dihedral angles of $\pi/3$.
Therefore, $\langle G_T,e_1,e_g\rangle=\Gamma^*$ is discrete
\cite{DM} and a quarter of $T$ is a fundamental polyhedron for
$\Gamma^*$. This is one of the groups from item~(2)({\it ii}) of the
theorem.

\medskip
In {\bf Case (5)}, $n\geq 5$, $m=n$, $l=n/4$. Analogously to
Case (1), there exists a reflection plane $\eta$ passing through
$e_1$ and $f$.
Conjugating $\eta$ by $R_\eta e_g$, we have a new reflection plane
through $V_1$.
However, the link of $V_1$ cannot be
divided by additional reflection planes so that the corresponding
reflection group to be discrete (see \cite{KS}).
Thus, $\Gamma^*$ is not discrete.

\medskip
In {\bf Case (6)}, $n=3$, $m=5$, $l=3/2$. Construct the decomposition
of $\widetilde\P$ into tetrahedra (see Figure~\ref{ptf_discr}d).
Note that $e_1$ coincides with one of the edges of a tetrahedron
in the decomposition. Consider the tetrahedron
$ABV_1V_2=[5,3/2,5;5,3/2,5]$. It is easy to see that the line
passing through the midpoints of $AV_1$ and $BV_2$ is orthogonal to
$\omega$ and to $AV_1$. Thus, it coincides with $e_g$. Moreover,
$ABV_1V_2$ is divided by reflection planes into four tetrahedra
$T=[2,3,5;2,3,2]$ so that $e_g$ passes through the midpoints of
the opposite edges of $T$ with dihedral angles of $\pi/5$ and $\pi/2$.
By \cite{DM}, $\langle G_T,e_g\rangle=\Gamma^*$ is discrete and a
half of $T$ is its fundamental polyhedron. This is one of the
groups from item~(2)({\it ii}) of the theorem.

\medskip
In {\bf Case (7)}, $n=5$, $m=3$, $l=3/2$. The bisector $\zeta$ is a
reflection plane and makes dihedral angles of $2\pi/5$ and
$3\pi/5$ with $\alpha$, and $\omega$, respectively. Hence, the link of
the
vertex $V_2$ formed by $\alpha$, $\omega$ and $\zeta$ is a spherical
triangle $(5,5/2,5/3)$. But the reflections in its sides do not
generate a discrete group~\cite{F2}.
Thus, $\Gamma^*$ is not discrete.

\medskip
In {\bf Case (8)}, $n=5$, $m=2$, $l=5/3$. In this case
$\widetilde\P$ is decomposed into tetrahedra $T=[2,3,5;2,3,2]$
and $e_1$ passes through the midpoints of the opposite edges of $T$
with dihedral angles of $\pi/5$ and $\pi/2$. One can see that
$e_g$ is conjugate to $e_1$ (see Figure~\ref{ptf_discr}e, where we
give the most important elements of the decomposition). From
\cite{DM}, $\langle G_T,e_g\rangle=\langle
G_T,e_1\rangle=\Gamma^*$ is discrete and a half of $T$ is its
fundamental polyhedron. This is one of the groups from item~(2)({\it
ii}) of the theorem.

\medskip
In {\bf Case (9)}, $n=3$, $m=5$, $l=5/4$, $\widetilde\P$ is
decomposed into tetrahedra $T=[2,3,5;2,2,4]$ (see
Figure~\ref{ptf_discr}f). Note that $e_1$ coincides with an edge of
$T$ and therefore, $e_1\in G_T$. Moreover, the edge $AV_1$ is
orthogonal to the edge $\alpha\cap\alpha'$ and lies in a plane
orthogonal to $\omega$. Hence, $e_g$ passes through $AV_1$ and
therefore, $e_g\in G_T$. Thus, $\Gamma^*=G_T$ is discrete; this is
one of the groups from item~(2)({\it ii}) of the theorem.

\medskip
In {\bf Case (10)}, $n=5$, $m=3$, $l=5/4$. In this case
$\widetilde\P$ is decomposed into tetrahedra $T=[2,3,5;2,2,5]$
(see Figure~\ref{ptf_discr}g). Since the bisector $\zeta$ is a
reflection plane, $e_1$ coincides with one of the edges of $T$ and
therefore, $e_1\in G_T$. Furthermore, $e_g$ passes through the
midpoints of the opposite edges of $T$ with dihedral angles of
$\pi/3$ and $\pi/2$. From \cite{DM}, $\Gamma^*=\langle
G_T,e_g\rangle$ is discrete and a half of $T$ is a fundamental
polyhedron for $\Gamma^*$. This is one of the groups from
item~(2)({\it ii}) of the theorem.

\medskip
In {\bf Case (11)}, $n=5$, $m=2$, $l=3/2$. In this case
$\widetilde\P$ is decomposed into tetrahedra $T=[2,3,5;2,3,2]$
(see Figure~\ref{ptf_discr}h). Again $e_1\in G_T$. It is not
difficult to see that $e_g$ passes through the midpoints of the
edges of $T$ with dihedral angles of $\pi/5$ and $\pi/2$.
$\Gamma^*=\langle G_T,e_g\rangle$ is discrete (see \cite{DM}) and
a half of $T$ is a fundamental polyhedron for $\Gamma^*$. This is
one of the groups from item~(2)({\it ii}) of the theorem.

\medskip
In {\bf Cases (12)} and {\bf (13)}, the link of a vertex formed by $\alpha$,
$\omega$, and $\zeta$ is a spherical triangle $(n,7/3,7/4)$.
Obviously, the group generated by reflections in sides of such a
triangle is not discrete. Thus, $\Gamma^*$ is not discrete.

\section{$m\in\{\infty,\overline\infty\}$ and $l$ is fractional $(l>1)$}

Suppose $m=\infty$.
Then the link of a vertex formed by $\alpha$, $\alpha'$,
and $\alpha''$ is a hyperbolic triangle $(n,l,\infty)$, where $l$ is
fractional. Clearly, there are no such discrete groups (cf.~\cite{Ma}).

\begin{figure}[htbp]
\centering
\begin{tabular}{cc}
\includegraphics[width=3.8 cm]{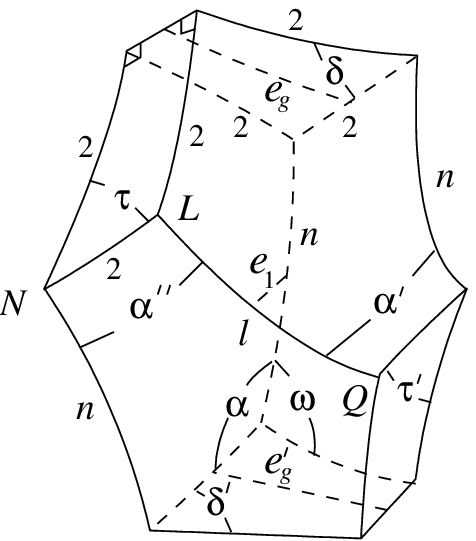}\qquad\quad &
\qquad \quad\includegraphics[width=2.6 cm]{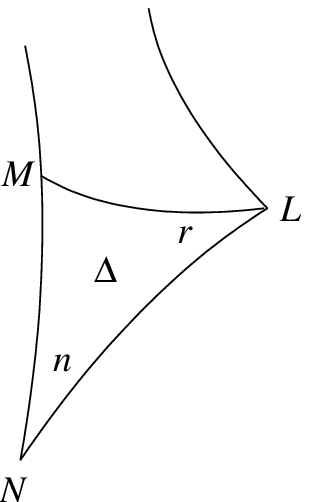}\\
(a)\qquad\quad &\qquad \quad (b)\\
\end{tabular}
\caption{}\label{polyhedron_pti}
\end{figure}

From here on, $m=\overline\infty$, $l$ is fractional, and $\Gamma^*$ is
discrete. Start with $\widetilde\P$ (see Section~5). The planes
$\alpha$, $\omega$, $\alpha'=e_g(\alpha)$, $\alpha''=e_1(\alpha')$
are reflection planes and make the following angles:
$\alpha$ and $\omega$, $\alpha$ and $\alpha''$, $\alpha'$ and
$\omega$ intersect at $\pi/n$, $\alpha'$ and $\alpha''$ at
$\pi/l$; but now the planes $\alpha$ and $\alpha'$ as well as $\alpha''$
and $\omega$ are disjoint. To draw this non-compact polyhedron
we use four additional planes to compactify it:
$\delta$ (which is orthogonal to $\alpha$, $\omega$, and $\alpha'$),
$\tau$ (orthogonal to $\alpha$, $\alpha'$, $\alpha''$),
$\delta'=e_1(\delta)$, and $\tau'=e_1(\tau)$ (see Figure~\ref{polyhedron_pti}a).

The compact polyhedron is also symmetric with respect to $e_1$.
Note that $\tau$ intersects $\delta$ and $p=\tau\cap\delta$
is the common perpendicular for $\alpha$ and $\alpha'$.
Remember that a priori $\tau$ and $\delta$ are not reflection
planes, so they can make any angle, which does not contradict
discreteness; moreover, $p$ can intersect $NL$ (as in
Case~(16) below).

Consider $\tau$ and
$\langle R_\alpha,R_\alpha',R_\alpha''\rangle$.
The group acts in $\tau$ as a group generated by reflections
in the sides of a polygon $D$ of infinite area with three sides and
two angles of $\pi/n$ and $\pi/l$ at vertices $N$ and $L$,
respectively. Since $l$ is fractional, there are additional
reflection lines in $\tau$ that correspond to reflection planes
through $LQ=\alpha'\cap\alpha''$ in $\widetilde\P$.

\medskip
\noindent
{\bf Remark~2.}
Note that if a plane through $LQ$ intersects $\alpha$, it also
intersects $\alpha\cap\tau$. Hence, if $\zeta$ is the bisector
of the dihedral angle of $\widetilde\P$ formed by $\alpha'$ and
$\alpha''$, it passes through $e_1$ and therefore,
 intersects $f$, $\alpha$,
and $\alpha\cap\tau$.

\medskip
\noindent
{\bf Remark~3.}
The bisector $\zeta$ from Remark~2
cannot be orthogonal to $\alpha$ (otherwise
$\alpha$ and $\alpha'$ meet at an angle of $\pi/n$, but they
must be disjoint).
If $\zeta$ is a reflection plane in $\Gamma^*$,
then $n=3$ or $5$ and the order of
$h=R_\alpha R_\zeta$ is $3$, $4$ or $5$ by Remark~1 (see Section~5).

\medskip
Now we will produce a complete list of possible decompositions
of $\tau$ by reflection planes.
Suppose $\pi/l=s\pi/r$, where $s$ and $r$ are
integers, $(s,r)=1$, $r>s\geq 2$. Then the line through $L$
cutting $D$ and making an angle of $\pi/r$ with $LN$ intersects
the opposite side $\alpha\cap\tau$ of $D$ at some point $M$
(because $\zeta$ also intersects $\alpha\cap\tau$ by Remark~2
and $s\geq 2$). Denote by $\Delta$ the triangle with vertices
$L$, $M$, $N$. It has two primitive angles $\pi/n$
and $\pi/r$ at $N$ and $L$, respectively, so either $\Delta$ is
one of the hyperbolic triangles in Figure~\ref{triangles}
or all its angles
are primitive. It is clear that the decomposition of $\tau$ is
specified by $n$, $l$, and $\Delta$.

Now suppose that $\Delta$ is chosen and placed into $\tau$ so that it
has primitive angles
$\phi=\pi/n$  ($n$ is odd) at $N$ and $\psi=\pi/r$ at $L$
(see Figure~\ref{polyhedron_pti}b). We find the
smallest $s$ such that $\alpha'$ and $\alpha$ are disjoint
and $(s,r)=1$. Then
we check if $\alpha$ and $\zeta$ intersect and if Remark~3 holds.
If both conditions hold,
we include such a decomposition in our list and replace
$s$ with $s+1$.
If just $\alpha\cap\zeta\neq\emptyset$, we do not include
the decomposition of $\tau$ in the list but still replace
$s$ with $s+1$.
Then check again if $(s,r)=1$, $\alpha\cap\zeta\neq\emptyset$,
and Remark~3 holds.
Using Remark~2, we stop our process for the fixed position of
$\Delta$ as soon as the planes $\alpha$ and $\zeta$ are parallel
or disjoint for some $s=\widetilde s$. We do not include
decompositions of $\tau$
with such a triangle for all $s\geq \widetilde s$.
Note that we can detect the behavior of
$\alpha$ and $\zeta$ (as well as $\alpha'$ and $\alpha$)
by considering their projections onto $\tau$, since $\alpha$, $\alpha'$,
and $\zeta$ are orthogonal to $\tau$.

Note that one picture in Figure~\ref{triangles} can give us
two different sequences of triangles, say for (h3), we can
take first $n=3$ and then $n=p\geq 7$.
If the same triangle can be placed in $\tau$ in a different way,
we proceed with a new position (in fact, we consider this new
position as a new triangle $\Delta$ in our algorithm).

It remains to explain why not all triangles with two primitive angles
give at least one decomposition of $\tau$ and appear as
$\Delta$ in our list.

We start with triangles (h2)--(h5) in Figure~\ref{triangles};
triangles of type (h1) will be considered later.
Obviously, triangles of type (h2) do not give any suitable
decomposition.
If a picture, say (h3), gives us a sequence of triangles that depends
on a parameter $p$, then we start with the smallest
$p=p_0$ (in fact, $p_0=7$) and consider this triangle as
$\Delta=\Delta(p_0)$. If for the chosen triangle the process
of producing the list of decompositions stops with some
$s=\widetilde s$, there is no need to consider all other
triangles
$\Delta(p)$
(of the same type and placed in a similar way
into $\tau$) for all $p\geq p_0$ and
$s\geq \widetilde s$ by Lemma~\ref{lemma2}.

Using the above algorithm, one can see that
if $\Delta$ is chosen from (h2)--(h5) in Figure~\ref{triangles},
then the only possible decompositions that satisfy both Remark~2
and the condition $\alpha\cap\zeta\neq\emptyset$
are those listed as (1)--(4) in Proposition~\ref{minfty_cases} below.
To produce the remaining cases in the list we need also the following

\begin{lemma}\label{lemma3}
If $\Delta$ has angles $\pi/5$ at $N$, $\pi/4$ at $L$, and $\pi/3$ at $M$,
then $\alpha\cap\zeta=\emptyset$ for $s=3$.
\end{lemma}

\noindent
{\it Proof.} The proof is a straightforward calculation.\qed

Now consider the two-parameter family of triangles (h1)
in Figure~\ref{triangles}.
Here $q$ and $p=n$ are odd $(1/q+1/p<1/2)$.
Fix $q=3$ first. Then we have a one-parameter family and the process
above gives decompositions (5) and (6) in
Proposition~\ref{minfty_cases}. For $q=5$, we obtain (7)
and (8).
For $q\geq 7$, using Remark~3 (when $s=2$) and
Lemmas~\ref{lemma2} and~\ref{lemma3} (when $s=3$), one can show that
the triangles (h1) do not give suitable decompositions.

It remains only to consider triangles with all primitive angles.
We start with a right-angled triangle $\Delta$, then $\angle M=\pi/2$.
If $n=3$ we use the above algorithm to get (9) and (10)
in the proposition below. Analogously,
using Remark~3 and Lemma~\ref{lemma2} for $n\geq 5$, we get (11) and
(12) and show that there are no other decompositions with right-angled
triangles.

Suppose that the greatest angle of $\Delta$ is $\pi/3$. It is easier to
start now with $s=3$ to see what happens. Since the angle at $L$ should be
less than $\pi$, i.e., $s\pi/r<\pi$, we have $r\geq 4$. Then
Lemma~\ref{lemma3} combined with Lemma~\ref{lemma2} gives that
we do not have other decompositions except (13) for $s=3$, and
therefore, for all $s\geq 3$. The last case to consider is $s=2$, that is the
case where the bisector $\zeta$ is the only reflection plane through
$LQ$. Using Remark~3 (which is now a criterion to stop),
we get the rest of the
list. We have proved the following

\begin{proposition}\label{minfty_cases}
Let $f$ and $g$ be as in Theorem~A.
If $\Gamma=\left<f,g\right>$ is discrete and $gfg^{-1}f$ is hyperbolic, then
either $l$ is an integer, or $n$, $l$, and $\Delta$ as above belong
to the following list:
\begin{itemize}
 \item[$(1)$] $n=3$, $l=8/3$, $\Delta=(3,8,8/3)$;
 \item[$(2)$] $n=3$, $l=7/3$, $\Delta=(3,7,7/3)$;
 \item[$(3)$] $n=7$, $l=7/3$, $\Delta=(7,7,7/4)$;
 \item[$(4)$] $n=7$, $l=7/5$, $\Delta=(7,7,7/4)$;
 \item[$(5)$] $n\geq 7$, $l=n/5$, $(n,5)=1$, $\Delta=(n,n,3/2)$;
 \item[$(6)$] $n\geq 9$, $l=n/7$, $(n,7)=1$, $\Delta=(n,n,3/2)$;
 \item[$(7)$] $n=5$, $l=5/2$, $\Delta=(5,5,5/2)$;
 \item[$(8)$] $n\geq 5$, $l=n/3$, $(n,3)=1$, $\Delta=(n,n,5/2)$;
 \item[$(9)$] $n=3$, $l=r/4$, $r\geq 7$, $(r,2)=1$, $\Delta=(3,r,2)$;
 \item[$(10)$] $n=3$, $l=r/5$, $r\geq 7$, $(r,5)=1$, $\Delta=(3,r,2)$;
 \item[$(11)$] $n\geq 5$, $l=r/3$, $r\geq 4$, $(r,3)=1$, $\Delta=(n,r,2)$;
 \item[$(12)$] $n=5$, $l=r/4$, $r\geq 5$, $(r,2)=1$, $\Delta=(5,r,2)$;
 \item[$(13)$] $n=3$, $l=r/3$, $r\geq 4$, $(r,3)=1$, $\Delta=(3,r,3)$;
 \item[$(14)$] $n=3$, $l=r/2$, $r\geq 5$, $(r,2)=1$, $\Delta=(3,r,3)$;
 \item[$(15)$] $n=3$, $l=r/2$, $r\geq 3$, $(r,2)=1$, $\Delta=(3,r,4)$;
 \item[$(16)$] $n=3$, $l=r/2$, $r\geq 3$, $(r,2)=1$, $\Delta=(3,r,5)$;
 \item[$(17)$] $n=5$, $l=r/2$, $r\geq 3$, $(r,2)=1$, $\Delta=(5,r,3)$;
 \item[$(18)$] $n=5$, $l=r/2$, $r\geq 3$, $(r,2)=1$, $\Delta=(5,r,4)$;
 \item[$(19)$] $n=5$, $l=r/2$, $r\geq 3$, $(r,2)=1$, $\Delta=(5,r,5)$.
\end{itemize}
\end{proposition}

Later we will see that Cases (16), (9), and (14)
(maybe with some further restrictions) lead to Item~({\it iii})
and Cases (8), (11) lead to Item~({\it iv})
of the theorem and that there are no other discrete groups among
(1)--(19).

\begin{figure}[htbp]
\centering
\begin{tabular}{cc}
\includegraphics[width=5 cm]{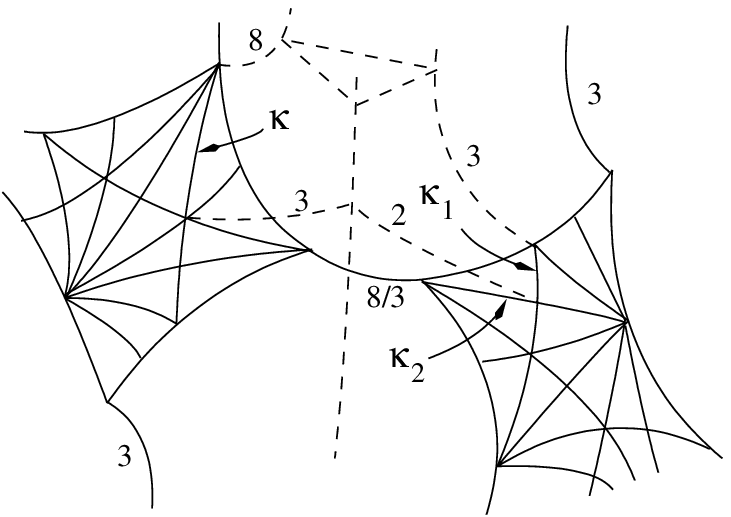}\qquad
& \qquad\includegraphics[width=4 cm]{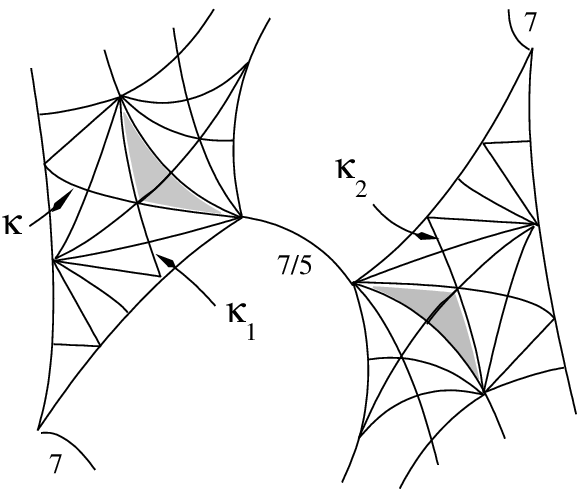}\\
(a)\qquad & \qquad(b)\\
\end{tabular}
 \caption{}\label{pti_1_4}
\end{figure}

\medskip
In {\bf Case (1)}, $n=3$, $l=8/3$, $\Delta=(3,8,8/3)$.
Let $\kappa$, $\kappa_1$, and $\kappa_2$ be three
planes in the decomposition of $\widetilde\P$ shown in
Figure~\ref{pti_1_4}a (if a plane, say $\kappa$, is orthogonal to $\tau$
or $\tau'=e_1(\tau)$, we show often just $\kappa\cap\tau$ and
$\kappa\cap\tau'$ in figures).

Suppose $\kappa$ and $\kappa_1$ are parallel or disjoint. Then
$\kappa_1$ is disjoint from $\alpha$ (because $\kappa$ and $\alpha$
are also disjoint and $\alpha$ and $\kappa_1$ lie in different
half-spaces bounded by $\kappa$). The same argument shows that
$\alpha$ and $\omega$ are disjoint (both $\alpha$ and $\omega$ are
disjoint from $\kappa_1$ and lie in different half-spaces with
respect to $\kappa_1$). This contradicts the fact that
$f=R_\omega R_\alpha$ is elliptic.

Suppose $\kappa$ and $\kappa_1$ intersect and form a dihedral angle of
$\pi/k$. We see that the link formed by $\alpha'$,
$\kappa$, and $\kappa_1$ is $(8,3,k)$. So, $k=8/3$
(Figure~\ref{triangles}(h3)) or $k$ is an integer, $k\geq 2$. The
case $k=8/3$ is impossible, because the link formed by $\kappa$,
$\kappa_1$, and $\kappa_2$ is a spherical triangle $(3,2,8/3)$
that cannot appear in a discrete group.

The case $k\geq 3$ is also
impossible, because two planes $\alpha$ and $\omega$ are disjoint
from $\kappa$ and $\kappa_1$, respectively, so even if $\alpha$
and $\omega$ intersect, the angle between them is smaller
than between $\kappa$ and $\kappa_1$.
If $k=2$ ($\kappa$ is orthogonal to $\kappa_1$), then it is easy
to see that $R_{\kappa}(\alpha)$ is disjoint from $\kappa_1$, so
$\alpha$ itself is also disjoint from $\kappa_1$. We see again
that $\alpha$ and $\omega$ are disjoint, because they lie in
different half-spaces with respect to $\kappa_1$ and $\kappa_1$
and $\omega$ are disjoint.

\medskip
In {\bf Case (4)}, $n=7$, $l=7/5$, $\Delta=(7,7,7/4)$. Consider $\kappa_1$ and
$\kappa_2=e_1(\kappa_1)$ (see Figure~\ref{pti_1_4}b). These two planes
are either disjoint or parallel (and then $\alpha$ and $\omega$ are
disjoint), or intersect at an angle of
$\pi/7$ to make the only possible compact link $(7,7,3/2)$
together with $\kappa$. But in the latter case
$\alpha$ and $\omega$ cannot make the same angle $\pi/7$
contradicting the fact that $n=7$.
Therefore, $\Gamma^*$ is not discrete.

\begin{figure}[htbp]
\centering
\begin{tabular}{cc}
\includegraphics[width=4 cm]{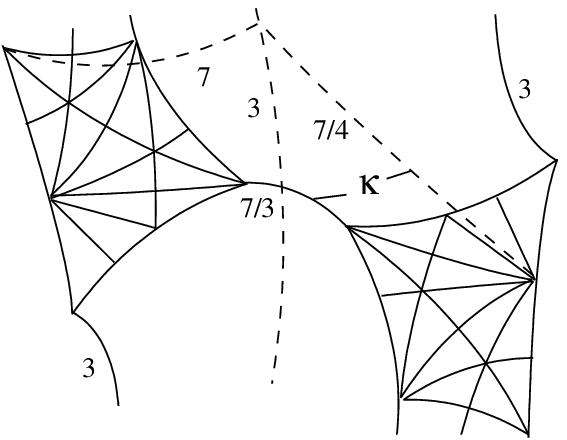} \qquad &
\qquad\includegraphics[width=4 cm]{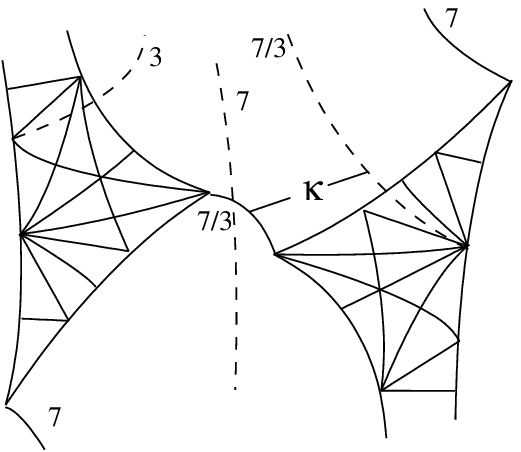}\\
(a)\qquad &\qquad (b)\\
\end{tabular}
\caption{}\label{pti_2_3}
\end{figure}

\medskip
In {\bf Cases (2)} and {\bf (3)}, $n=3$, $\Delta=(3,7,7/3)$
or $n=7$, $\Delta=(7,7,7/4)$ and $l=7/3$. Let $\kappa$ be a
plane in the decomposition of $\widetilde\P$ making a dihedral
angle of $\pi/7$ with $\alpha'$ (see Figure~\ref{pti_2_3}).  The link
formed by $\kappa$, $\alpha$, and $\omega$ is a
spherical triangle $(3,7,7/4)$
for Case (2) or a hyperbolic triangle $(3,7,7/3)$
for Case (3). However, the group generated by
reflections in its sides is not discrete (there are no such
triangles in Figure~\ref{triangles}).
Thus, in both cases $\Gamma^*$ is not
discrete.

\begin{figure}[htbp]
\centering
\includegraphics[width=4.6 cm]{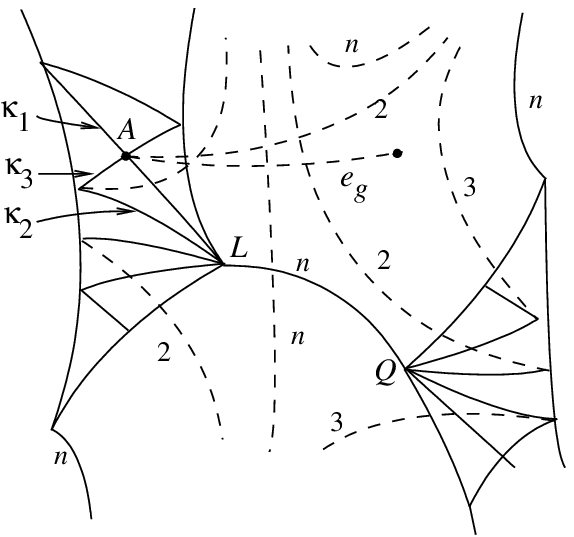}
\caption{}\label{pti_5}
\end{figure}

\medskip
In {\bf Case (5)}, $n\geq 7$, $l=n/5$, $(n,5)=1$, $\Delta=(n,n,3/2)$;
$\widetilde\P$ is
decomposed into infinite volume tetrahedra $T=[2,3,n;2,3,n]$ (see
Figure~\ref{pti_5}).

Now we determine the position of the axis $e_g$. It passes through
$A$ and is orthogonal to $\omega$. Since $e_g$ is not orthogonal
to $\tau$ (the perpendicular to $\tau$ through $A$ is the line
$\kappa_1\cap\kappa_3$), the distance between $e_g$ and
the axis $LQ=\alpha'\cap\alpha''$ of an elliptic element of order
$n$ is less than $AL$, where
$$\cosh AL=1/(2\sin(\pi/n))=\cosh\rho_{min}(2,n),\quad n\geq 7.$$
We arrive at a contradiction; $\Gamma^*$ is not discrete.

\begin{figure}[htbp]
\centering
\begin{tabular}{cc}
\includegraphics[width=4 cm]{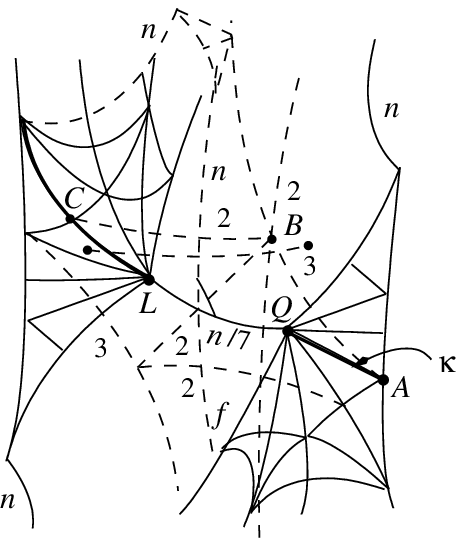}\qquad &
\qquad \includegraphics[width=3.5 cm]{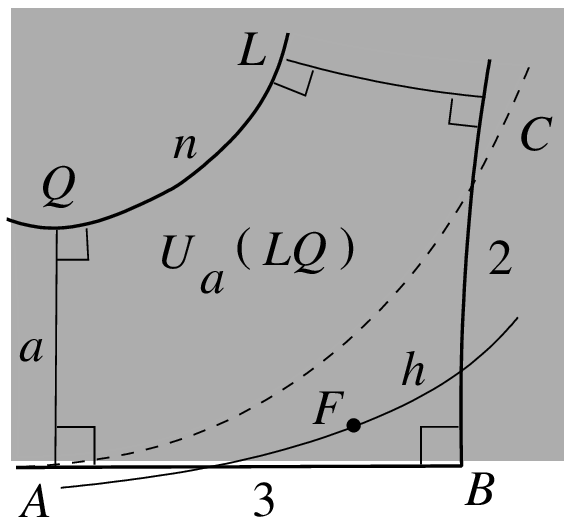}\\
(a)\qquad & \qquad (b)\\
\end{tabular}
\caption{}\label{pti_6}
\end{figure}

\medskip
In {\bf Case (6)}, $n\geq 9$, $l=n/7$, $(n,7)=1$, $\Delta=(n,n,3/2)$.
Construct a decomposition of
$\widetilde\P$ into infinite volume tetrahedra $T=[2,2,3;2,n,3]$ (see
Figure~\ref{pti_6}a).

The points of intersection of $e_g$ with $\tau$ and $\omega$ lie
in different half-spaces bounded by the plane $\kappa$. So, $e_g$
intersects $\kappa$. Moreover, the point $F=e_g\cap\kappa$ lies within
the pentagon $ABCLQ$ in Figure~\ref{pti_6}b. Since $e_g$ makes an
angle of less than $\pi/6$ with $\kappa$, there exists an elliptic
element $h\in \langle G_T,e_g\rangle$ of order $q>3$ with the axis
lying in $\kappa$.

On one hand, $h$ cannot intersect the $a$-neighborhood $U_a(LQ)$ of
the axis of order $n$ passing through $LQ$, where
$a=AQ=\rho_{min}(3,n)<\rho_{min}(q,n)$, for all $q>3$ (see
Table~\ref{table_min} and the formulas~(\ref{mindist}) or
\cite{GMM}). On the other hand, $h$ cannot meet $AB$
since the projection of $e_g$ on $\kappa$ is orthogonal to $h$ and
$AB$. Contradiction. $\Gamma^*$ is not discrete.

\medskip
In {\bf Case (7)}, $n=5$, $l=5/2$, $\Delta=(5,5,5/2)$.
The link of the vertex formed by $\alpha$, $\omega$, and $\zeta$,
where $\zeta$ is the bisector of the dihedral angle formed by
$\alpha'$ with $\alpha''$, is a spherical triangle $(5,5/2,5/3)$. But
the group generated by reflections in its sides is not discrete
\cite{F2}. Thus, $\Gamma^*$ is not discrete.

\begin{figure}[htbp]
\centering
\begin{tabular}{cc}
\includegraphics[width=5.5 cm]{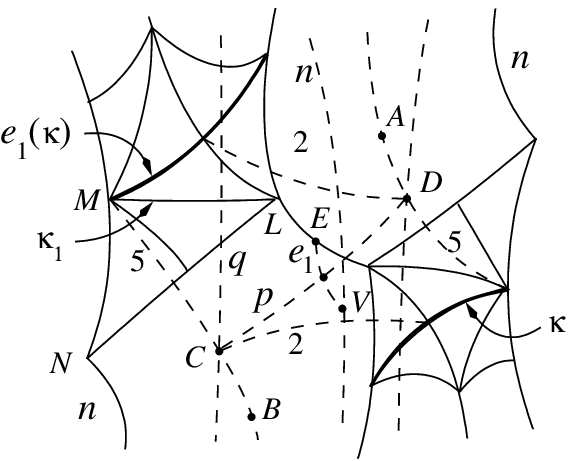}\qquad &
\qquad \includegraphics[width=5 cm]{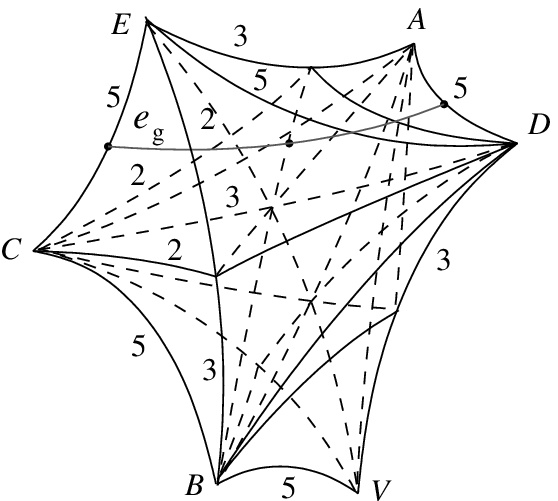}
\\
(a)\qquad & \qquad (b)\\
\end{tabular}
\caption{}\label{pti_8}
\end{figure}

\medskip
In {\bf Case (8)}, $n\geq 5$, $l=n/3$, $(n,3)=1$, $\Delta=(n,n,5/2)$.
Consider two
reflection planes $\kappa$ and $e_1(\kappa)$ from the decomposition
(see Figure~\ref{pti_8}a).
Denote by $\pi/p$ and $\pi/q$ the angles that $\kappa$ makes with
$e_1(\kappa)$ and $\alpha$, respectively.
It is not difficult to show that if $\pi/p$
is less than or equal to $\pi/2$ (or the
planes are parallel or disjoint) then $\alpha$ and $\omega$ are
disjoint.

So we can assume that $\pi/p>\pi/2$.
The link of the vertex formed
by $\alpha$, $\omega$, and $\kappa$ is a triangle
$\Delta=(5/2,n,q)$.
Suppose $\Delta$
is spherical, then immediately $n=5$ and from Figure~\ref{triangles},
$q=2$ or $q=3/2$. Therefore, from the link
of $C$, either $p=2$ for $q=2$ or $p=5/2$ for $q=3/2$.
In both cases $\pi/p\leq\pi/2$ and thus, $\alpha$ and $\omega$ are
disjoint, which is impossible.
Suppose now that $\Delta$ is hyperbolic. By~\cite[Lemma~2.1]{KS},
$q=n\geq 5$ and the link of
$C$ made by $\alpha$, $\kappa$, and $\kappa_1$ is a spherical triangle
$(2,5/3,q/(q-1))$. Then $q=n=5$ and therefore, $p=3/2$ (see \cite{F2}).

Consider $n=5$, $p=3/2$. Then $\widetilde\P$ can be decomposed
into compact tetrahedra
$T=[2,2,3;2,5,3]$ (see Figure~\ref{pti_8}b). Here, $e_1$ coincides
with an edge of $T$
(actually, with $EV$). Determine the position of $e_g$.
For $n=5$, $e_g$ intersects $LM$ in the midpoint.
Further, $CE$ makes equal alternate interior angles with
$\kappa_1\cap\alpha$ and $\kappa_1\cap\alpha'$ and hence,
$LM$ bisects $CE$. Therefore, $e_g$ also bisects $CE$.
The tetrahedron $ACDE=[5,5,3/2;5,5,3/2]$ consists of four copies of $T$.
Consider the line $l$ passing through the middles of $CE$ and $AD$
in $ACDE$. Then $l$
intersects an edge with the dihedral angle $\pi/2$ that is opposite
to both $CE$ and $AD$. Thus, $l$ is orthogonal to $CE$ and $AD$.
Moreover, $l$ bisects the dihedral angle formed by $\kappa=ACD$ and
$e_1(\kappa_1)=AED$ and therefore, $l$ is orthogonal to $\omega$.
So, $e_g$ coincides with $l$.
By \cite{DM}, $\Gamma^*=\langle G_T,e_g\rangle$ is discrete.
This is one of the groups from
item~(2)({\it iv}) of the theorem.

\begin{figure}[htbp]
\centering
\includegraphics[width=4 cm]{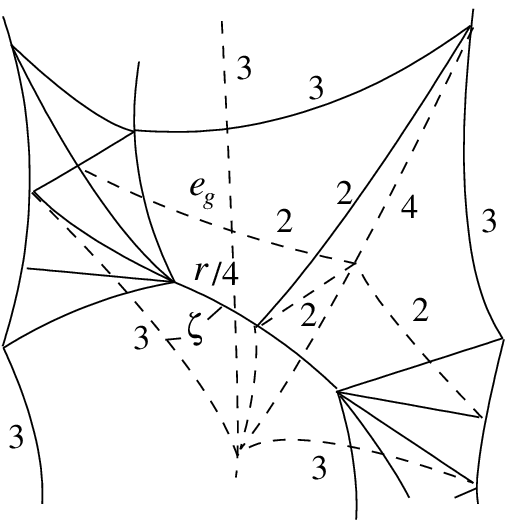}
\caption{}\label{pti_9}
\end{figure}

\medskip
In {\bf Case (9)}, $n=3$, $l=r/4$, $r\geq 7$, $(r,2)=1$,
$\Delta=(3,r,2)$.
Taking into
account the fact that the bisector $\zeta$ is a reflection plane,
it is not difficult to construct the decomposition of $\widetilde\P$
into tetrahedra $T=[2,3,r;2,2,4]$  with one imaginary
vertex. By construction, $e_1$ coincides with an edge of $T$ with
the dihedral angle of $\pi/2$. Moreover, one can see that $e_g$
coincides with another edge of $T$ with the dihedral angle of
$\pi/2$ (see Figure~\ref{pti_9}). Hence, $\Gamma^*=G_T$
is discrete and $T$ is its fundamental
polyhedron. It is one of the groups from item~(2)({\it iii}) of the
theorem.

\medskip
In {\bf Case (10)}, $n=3$, $l=r/5$, $r\geq 7$, $(r,5)=1$,
$\Delta=(3,r,2)$. There is a
plane $\kappa$ from the decomposition of $\widetilde\P$ which
passes through $\alpha'\cap\alpha''$ and makes an angle of
$2\pi/r$ with $\alpha'$. The link of the vertex formed by $\alpha$,
$\omega$, and $\kappa$ is a spherical triangle $(3,r,3/2)$, $r\geq
7$. But the group generated by reflections in the sides of this
triangle is not discrete. Thus, $\Gamma^*$ is not discrete.

\begin{figure}[htbp]
\centering
\begin{tabular}{cc}
\includegraphics[width=3.2 cm]{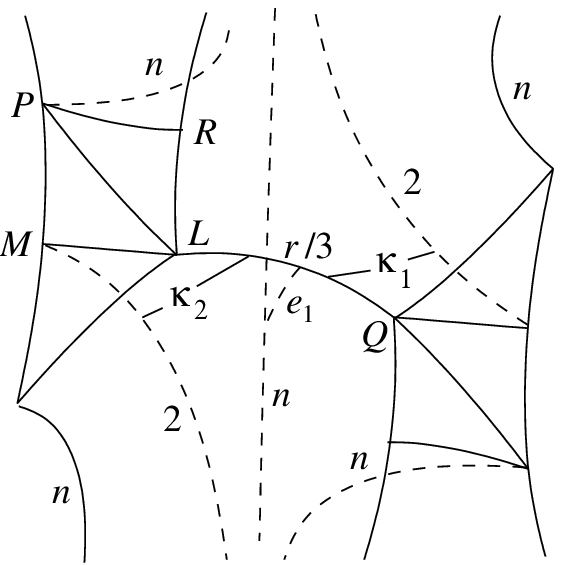}\qquad
& \qquad\includegraphics[width=4 cm]{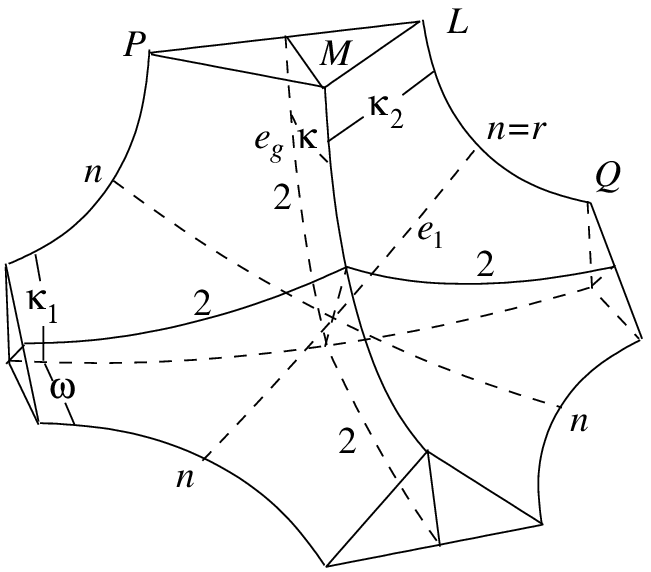}\\
(a)\qquad &\qquad (b)\\
\end{tabular}
\caption{}\label{pti_11}
\end{figure}

\medskip
In {\bf Case (11)}, $n\geq 5$, $l=r/3$, $r\geq 4$, $(r,3)=1$,
$\Delta=(n,r,2)$. One can
construct the decomposition of $\widetilde\P$ into tetrahedra
$T=[2,n,n,;2,n,r]$ of infinite volume (see Figure~\ref{pti_11}).
Let $\kappa_1$ and $\kappa_2$ be reflection planes that pass through
$\alpha'\cap\alpha''$ and make angles of $\pi/r$ and $2\pi/r$ with
$\alpha'$, respectively. Then $T$ is bounded by $\alpha$,
$\kappa_1$, $\kappa_2$, and $\omega$.

\smallskip
\noindent
(1) Assume $n=r\geq 5$. Then $e_g$ lies in $\kappa_1$. Moreover, $e_g$
maps the face of $T$ lying in $\kappa_1$ onto itself. There are two
more reflection planes that meet the interior of $T$. One of them,
denote it by $\kappa$,
passes through $e_g$ and is orthogonal to $\kappa_1$ ($\kappa$
passes also through $\alpha\cap\kappa_2$), its image with respect to
$e_1$ is the other reflection plane that passes through
$\kappa_1\cap\omega$.  Therefore, $T$ is decomposed into four
tetrahedra $T_1=[2,2,4;2,n,4]$ (see Figure~\ref{pti_11}b).
Obviously, $\Gamma^*=\langle G_{T_1},e_1\rangle$ is discrete and a
half of $T_1$ is a fundamental polyhedron for $\Gamma^*$. This is
item~(2)$(iv)$ of Theorem~A.

\begin{figure}[htbp]
\centering
\includegraphics[width=4 cm]{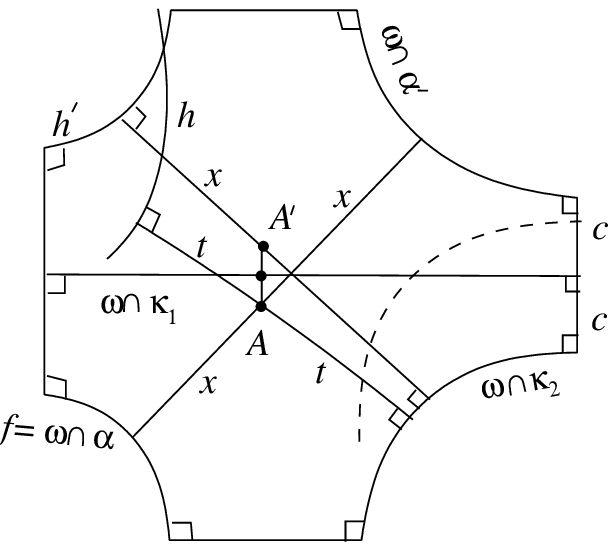}
\caption{}\label{pti_11_omega}
\end{figure}

\smallskip
\noindent
(2) Assume $n>r\geq 4$.
In this case $e_g$ does not lie in a
reflection plane. Then $e_g$ maps $\kappa_2$ to a plane $\kappa_3$;
moreover,
$\kappa_3$ makes an angle of $\pi/n$ with $\omega$. Let
$h=e_gR_\omega R_{\kappa_2}e_g$ and
let $e'_g=R_{\kappa_1}e_gR_{\kappa_1}$.
Let $A=e_g\cap\omega$ and $A'=e'_g\cap\omega$.
The distance from $A$ to $h$ (denote it by $t$)
is less than the
distance from $A'$ to the axis of
$h'=e'_gR_\omega R_{\kappa_2}e'_g$ (denote it by $x$,
see Figure~\ref{pti_11_omega}).

If $h$ and $h'$ intersect, then for $n\geq 7$, $\Gamma^*$
is not discrete. Suppose $n=5$ and $h$ and $h'$ intersect at an angle
$\phi$.

Elementary calculations show that for $n=5$, $r=4$, we have
$1>\cos\phi>\cos(\pi/6)$, that is $\Gamma^*$ is not discrete.

So we may assume that $n>5$ and $h$ and $h'$ are disjoint.
Let $h_1=e'_ghe'_g$. The elliptic element $h_1$ has order $n>5$.
The axis $h_1$ then cannot intersect (or be parallel to)
the axes of order $n$ through $\omega\cap\kappa_2$ and
$\omega\cap\alpha'$ ($h_1$ is drawn in Figure~\ref{pti_11_omega}
by the dashed line). Thus,
$$
\rho(h_1,\omega\cap\kappa_2)+\rho(h_1,\omega\cap\alpha')<2c,
$$
where $c$ can be defined from the link of the imaginary vertex made by
$\kappa_1$, $\kappa_2$, and $\omega$.
However,
$$
c={\rm arccosh}\left(\frac{\cos(\pi/r)}{\sin(\pi/n)}\right)<
{\rm arccosh}\left(\frac{\cos(2\pi/n)}{2\sin^2(\pi/n)}\right)=
\rho_{min}(n,n)
$$
for all $n > r \geq 8$.
We arrive at a contradiction, $\Gamma^*$ is not discrete.

\smallskip
\noindent
(c) If $r>n\geq 5$ the arguments are analogous and the group
is also not discrete.

\begin{figure}[htbp]
\centering
\includegraphics[width=4 cm]{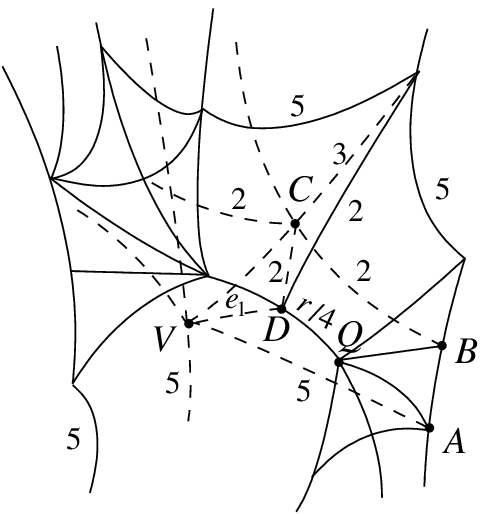}
\caption{}\label{pti_12}
\end{figure}

\medskip
In {\bf Case (12)}, $n=5$, $l=r/4$, $r\geq 5$, $(r,2)=1$,
$\Delta=(5,r,2)$. In this case
$\widetilde\P$ is decomposed into tetrahedra $T=[2,2,3;2,5,r]$
of infinite volume (see Figure~\ref{pti_12}). Note that $e_1$
coincides with the edge $DV$ of $T$.

Elementary calculations show that $e_g$ intersects $\omega$
orthogonally inside the quadrilateral $ABCV$. Besides, $e_g$
intersects either the triangular face $CDV$ or
the quadrilateral $AVDQ$. Show that
in either case $\Gamma^*$ is not discrete.

\smallskip
\noindent (a) The axis $e_g$ intersects $CDV$. Then the angle
of intersection is less than $\pi/6$ and hence, there exists an
elliptic element $h\in \Gamma^*$ of order $q>3$ with the axis
lying in $CDV$. Since $h$ and $CV$ are disjoint, $h$ intersects
$CD$ and $DV$.

Let $H=h\cap CD$. Then $CD=CH+HD$. From $\triangle CDV$,
$\cosh CD=2\cos(\pi/5)/(\sqrt{3}\sin(\pi/r))$. On the other hand, the
distance between $h$ and an elliptic element of order $3$ or $r\geq 5$
cannot be greater than the corresponding minimal distance between
axes of elliptic
elements in a discrete group and therefore,
\begin{equation}\label{star}
CH\geq \rho_{min}(q,3) \qquad {\rm and} \qquad
HD\geq \rho_{min}(q,r).
\end{equation}
However,
\begin{equation}\label{twostar}
\frac{2\cos(\pi/5)}{\sqrt{3}\sin(\pi/r)}<
\cosh(\rho_{min}(q,r)+\rho_{min}(q,3)), \qquad {\rm for\ }
q>3, r\geq 5.
\end{equation}
Therefore, from (\ref{star})--(\ref{twostar}) it follows that
$CD<CH+HD$. Contradiction. So, $\Gamma^*$ is not discrete.

\smallskip
\noindent (b) The axis $e_g$ intersects $AVDQ$. Then there
exists an elliptic element $h\in\Gamma^*$ of order $q$ with the
axis lying in $AVDQ$ such that $h$ does not intersect $AV$, and
there exists a reflection plane $\kappa$ that passes through
$h$ and cuts the interior of $ABCVDQ$.
Consider the following situations:

\smallskip
\noindent (i) if $h$ and $DQ$ are disjoint and $q\geq 3$, then we
estimate the distances between $h$ and the axes of orders $5$ and
$r$ as above and conclude that the group is not discrete.

\smallskip
\noindent (ii) if $h$ and $DQ$ are disjoint and $q=2$, then the
line of intersection of $\kappa$ and $BCDQ$ is the axis of an
elliptic element of order $\widetilde q\geq 3$.
Similarly to part~(a) above,
$$
\cosh AQ=\cot(\pi/5)\cot(\pi/r)<
\cosh(\rho_{min}(\widetilde q,5)+\rho_{min}(\widetilde q,r))
$$
for all $\widetilde q\geq 3$, $r\geq 5$. Therefore, the group is
not discrete.

\smallskip
\noindent (iii) if $h$ intersects $DQ$ at some point (possibly at
infinity), then $r=5$ and $\kappa$ cuts off a finite volume
tetrahedron from the trihedral angle with the vertex $D$. Since
$\Gamma^*$ is discrete, $\langle G_T, R_\kappa\rangle$ is also
discrete and ${\Bbb H}^3$ is tesselated by finite volume Coxeter
tetrahedra. There are nine compact Coxeter tetrahedra, but the group
generated by reflections in faces of any of them does not contain
subgroups generated by reflections in the sides of triangles $(2,2,3)$
and $(2,2,5)$ simultaneously. As for non-compact tetrahedra, only
$T[2,2,5;2,3,6]$ contains such subgroups, but in this case an axis
of order $6$ passes through $CV$ and thus, the stabilizer of $V$
contains elements of orders $5$ and $6$, which is impossible in a
discrete group. So, $\Gamma^*$ is not discrete.

\begin{figure}[htbp]
\centering
\includegraphics[width=4 cm]{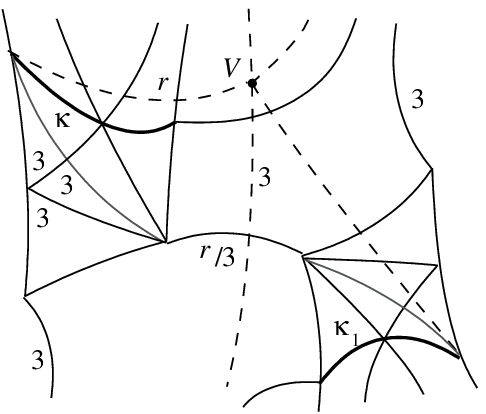}
\caption{}\label{pti_13}
\end{figure}

\medskip
In {\bf Case (13)}, $n=3$, $l=r/3$, $r\geq 4$, $(r,3)=1$,
$\Delta=(3,r,3)$. In this case
the bisector $\zeta$ is not a reflection plane. Let $V=\zeta\cap
f$ and let $\kappa$ and $\kappa_1$ be reflection planes such that
$\kappa$ passes through $\alpha\cap\zeta$ and makes an angle of
$\pi/r$ with $\alpha$ and $\kappa_1=e_1(\kappa)$
(Figure~\ref{pti_13}). Since $\zeta$ intersects $f$,
$\kappa$ also intersects $f$ at the same point $V$. Since the
elliptic element $R_\alpha R_{\kappa}$ of order $r$ belongs to
$St_{\Gamma^*}(V)\subset O(3)$, we have $r=4,5$.

Suppose $r=5$. Consider the link $\Delta$ of $V$ formed by
$\alpha$, $\omega$, and $\kappa$. Then $\Delta=(3,5/4,p)$, $p>0$.
Since the group is supposed to be discrete, $p=2,5/2,3$, or $5$.
But this means that the area of $\Delta$ is at least a half of the
area of the digonal link of $V$ formed by $\alpha$ and $\omega$,
which is impossible.
Arguments for $r=4$ are analogous and $\Gamma^*$ is not discrete.

\begin{figure}[htbp]
\centering
\includegraphics[width=3.5 cm]{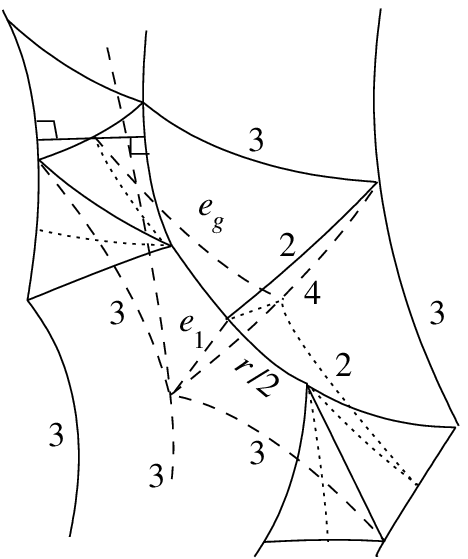}
\caption{}\label{pti_14}
\end{figure}

\medskip
In {\bf Case (14)}, $n=3$, $l=r/2$, $r\geq 5$, $(r,2)=1$,
$\Delta=(3,r,3)$. It is easy
to construct the decomposition of $\widetilde\P$ into tetrahedra
$T=[2,2,4;3,3,r]$ of infinite volume (see Figure~\ref{pti_14}).

In this decomposition, $e_1$ coincides with one of the edges of
$T$ and $e_g$ lies in the face $F$ of $T$ opposite to the vertex
$(2,2,r)$. Moreover, $e_g$ maps $F$ onto itself. Hence, there is
a reflection plane $\kappa$ passing through $e_g$ orthogonally to
$F$ so that $\kappa$ decomposes $T$ into two tetrahedra
$T_1=[2,2,4;2,3,2r]$ of infinite volume (see Figure~\ref{pti_14},
where new lines in the decomposition are dotted,
compare with Case (9)). Thus,
$\Gamma^*=G_{T_1}$ is discrete and is one of the groups from
item~(2)({\it iii}) of the theorem.

\begin{figure}[htbp]
\centering
\includegraphics[width=4.5 cm]{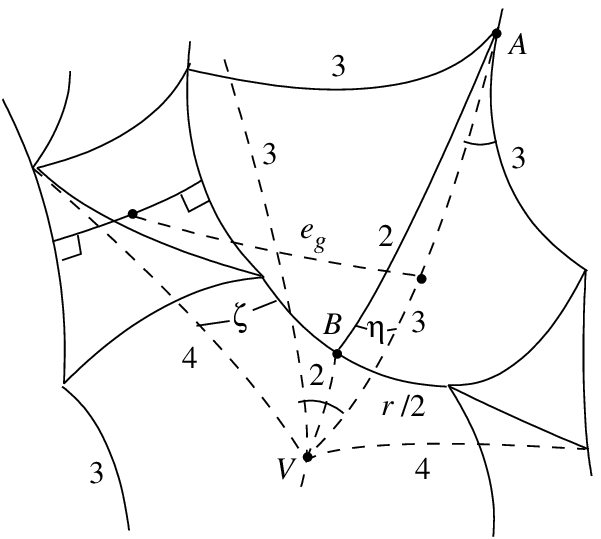}
\caption{}\label{pti_15}
\end{figure}

\medskip
In {\bf Case (15)}, $n=3$, $l=r/2$, $r\geq 3$, $(r,2)=1$,
$\Delta=(3,r,4)$. One can
construct the decomposition of $\widetilde\P$ into infinite volume
tetrahedra
$T=[2,2,3;3,4,r]$ (see Figure~\ref{pti_15}). In this decomposition,
$e_1$ coincides with one of the edges of $T$. We determine the
position of $e_g$. Let $\eta$ be the plane in the decomposition
that is orthogonal to the bisector $\zeta$. Let
$A=\eta\cap\alpha'\cap\omega$, $B=\eta\cap\alpha'\cap\zeta$, and
$V=\eta\cap\zeta\cap\omega$. Consider the face $F$ of
$\widetilde\P$ lying in $\omega$. $F$ is split by the axis $AV$ of
order $3$ into two parts so that $AV$ makes equal (alternate
interior) angles with $f=\alpha\cap\omega$ and
$\alpha'\cap\omega$. Then $e_g$ intersects $AV$ in the
midpoint.

Since $e_g$ is orthogonal to $AV$, it maps $A$ to $V$ and
therefore, the stabilizer $St_{\Gamma^*}(A)$
of the point $A$ in $\Gamma^*$ contains an element of order $4$.
Hence, $AB$
is the axis of an elliptic element of order $4$ and from the link
of $B$ formed by $\alpha'$, $\zeta$, and $\eta$ we see that $r=4$,
which contradicts $(r,2)=1$. Thus, $\Gamma^*$ is not discrete.

\begin{figure}[htbp]
\centering
\begin{tabular}{cc}
\includegraphics[width=4.5 cm]{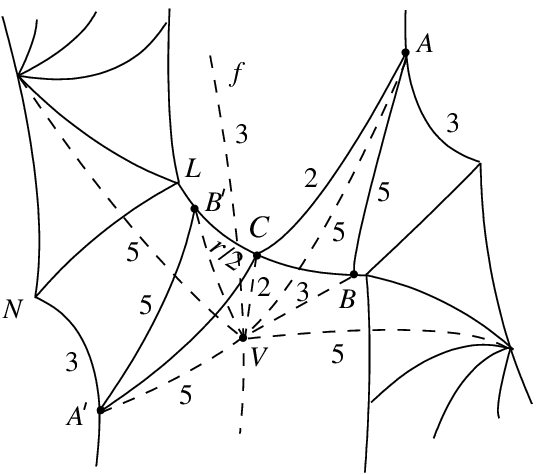}\quad &
\quad \includegraphics[width=4.5 cm]{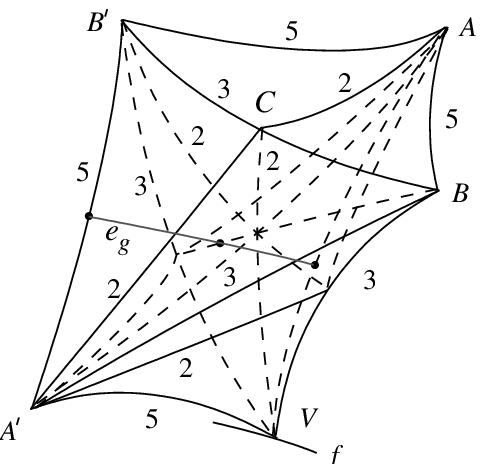}\\
(a)\quad & \quad (b)\\
\end{tabular}
\caption{}\label{pti_16}
\end{figure}

\medskip
In {\bf Case (16)}, $n=3$, $l=r/2$, $r\geq 3$, $(r,2)=1$,
$\Delta=(3,r,5)$.
Since the bisector $\zeta$ is a reflection plane, there is a
reflection plane $\eta$ passing through $e_1$ orthogonally to
$\zeta$ (see Figure~\ref{pti_16}a).
The link of the vertex $A$ formed by $\eta$, $\alpha'$,
and $\omega$ is a spherical triangle $(2,3,5/2)$. Hence, there is
a reflection plane $\kappa$ through $\eta\cap\omega$ that makes an
angle of $\pi/5$ with both $\omega$ and $\eta$. One can see that
the link of the vertex $B$ formed by $\alpha'$, $\alpha''$, and
$\kappa$ is a triangle $(3,r,5/2)$, $r\geq 3$. However, if $r\geq
5$, then $(3,r,5/2)$ is a hyperbolic triangle that does not
lead to a discrete group.

Consider $r=3$. Then the tetrahedron $ABCV$ is compact and
splits up into three tetrahedra $T=[2,2,3;2,5,3]$ (Figure~\ref{pti_16}b).
In the decomposition of $\widetilde\P$ into $T$, $e_1$ coincides
with an edge of $T$.

We determine
the position of $e_g$. First, note that for $r=3$, $e_g$ intersects
$LN$ in the midpoint. Further, $A'B'$ makes equal alternate interior
angles with $\alpha\cap\alpha'$ and $\alpha'\cap\alpha''$ and hence,
$LN$ bisects $A'B'$. Therefore, $e_g$ also bisects $A'B'$.

Consider the line $l$ passing through the midpoints of $A'B'$ and $AV$.
It bisects the edge with the dihedral angle $\pi/2$ that is opposite
to both $AV$ and $A'B'$ in the small tetrahedra. Thus, $l$ is
orthogonal to $AV$ and to $A'B'$. Moreover, $l$ makes an angle
of $\pi/10$ with $\eta=AA'V$ and, therefore, is orthogonal to
$\omega$. So, $e_g$ coincides with $l$.

By \cite{DM}, $\langle G_T,e_g\rangle=\Gamma^*$ is
discrete. This is one of the groups from item~(2)({\it iii}) of the
theorem.

\begin{figure}[htbp]
\centering
\begin{tabular}{cc}
\includegraphics[width=4.2 cm]{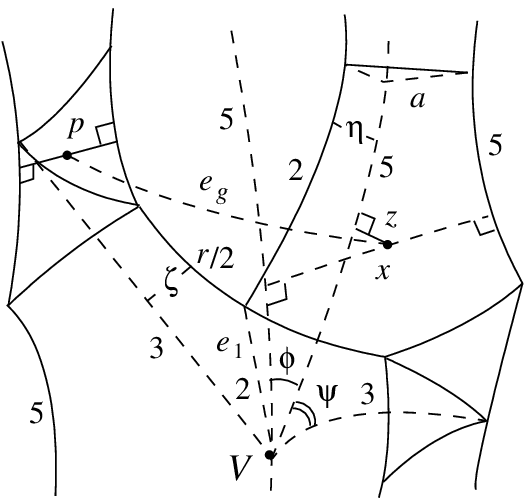} \qquad &
\qquad\includegraphics[width=4.2 cm]{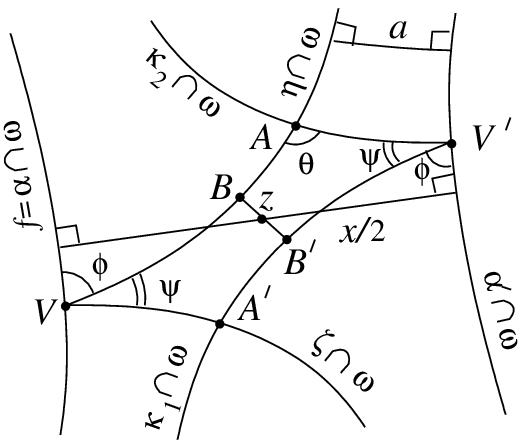}\\
(a)\qquad &\qquad (b)\\
\end{tabular}
\caption{}\label{pti_17}
\end{figure}

\medskip
In {\bf Case (17)}, $n=5$, $l=r/2$, $r\geq 3$, $(r,2)=1$,
$\Delta=(5,r,3)$. In this case
$\widetilde\P$ is decomposed into tetrahedra $T=[2,2,5;3,5,r]$
of infinite volume (see Figure~\ref{pti_17}a).

Let $z$ denote the distance between $e_g$ and $\eta\cap\omega$,
an axis of order~$5$.
From the decomposition of $\tau$, we can calculate the length of
the common perpendicular $p$ to $\alpha$ and $\alpha'$:
$$
\cosh p=\cos(\pi/r)+\cos(\pi/5).
$$
Then from the plane $\delta$, which is orthogonal to $\alpha$,
$\alpha'$, and $\omega$, we calculate the length of the common
perpendicular $x$ to $f=\alpha\cap\omega$ and $\alpha'\cap\omega$:
$$
\cosh x=\frac{\cosh p+\cos^2(\pi/5)}{\sin^2(\pi/5)}.
$$
From the link of the vertex $V$ formed by $\zeta$, $\alpha$, and
$\omega$,
$$
\cos \phi=\frac{\cos(2\pi/3)+\cos^2(\pi/5)}{\sin^2(\pi/5)}=
\frac{2\cos^2(\pi/5)-1}{2\sin^2(\pi/5)}.
$$
Analogously, from the link of a vertex formed by $\alpha'$,
$\eta$, $\omega$,
$$
\cosh a=\cot^2(\pi/5).
$$
Then
$$
\sinh z=\frac{\cosh a-\cos\phi}{2\sinh(x/2)}
=\frac{1}{4\sin^2(\pi/5)\sinh (x/2)}.
$$

If $r>3$, then $z<\rho_{min}(2,5)$ and the group is not discrete.
Assume that $r=3$. Let $\kappa_1$ and $\kappa_2$ be reflection planes such
that $R_{\kappa_1}=e_gR_\eta e_g$ and $R_{\kappa_2}=e_g R_\zeta e_g$.
Then $\eta$, $\kappa_1$, $\kappa_2$, and $\zeta$ cut off a quadrilateral
$V_1AV'_1A'$ from the plane $\omega$ (Figure~\ref{pti_17}b).
Calculate $\theta=\angle A$
in $VAV'A'$:
$$
\cos\theta=-\cos\psi\cosh(2z)+\sin\psi\sinh(2z)\sinh B'V',
$$
where
$$
\cos\psi=\cot(\pi/5)/\sqrt{3}\quad {\rm and} \quad
\sinh B'V'=\frac{\sinh(x/2)-\cos\phi\sinh z}{\sin\phi\cosh z}.
$$

On the other hand, the link of the vertex $A$ formed by $\eta$,
$\kappa_2$, and $\omega$ is a spherical triangle $(3,5,q)$. In a
discrete group, $q=2,5/3,3/2$, or $5/4$. Then from such a link
$$
\cos\theta(q)
=\frac{\cos(\pi/q)+\cos(\pi/3)\cos(\pi/5)}{\sin(\pi/3)\sin(\pi/5)}=
\frac{2\cos(\pi/q)+\cos(\pi/5)}{\sqrt{3}\sin(\pi/5)}
\neq\cos\theta
$$
for all possible $q$. We arrive at a contradiction and thus,
$\Gamma^*$ is not discrete.

\medskip
In {\bf Case (18)}, $n=5$, $l=r/2$, $r\geq 3$, $(r,2)=1$,
$\Delta=(5,r,4)$. The link of
a vertex formed by $\alpha$, $\omega$, and the bisector $\zeta$ is
a spherical triangle $(5,4,4/3)$. But the
group generated by reflections in its sides
is not discrete. Thus, $\Gamma^*$ is not discrete.

\begin{figure}[htbp]
\centering
\begin{tabular}{cc}
\includegraphics[width=4.5 cm]{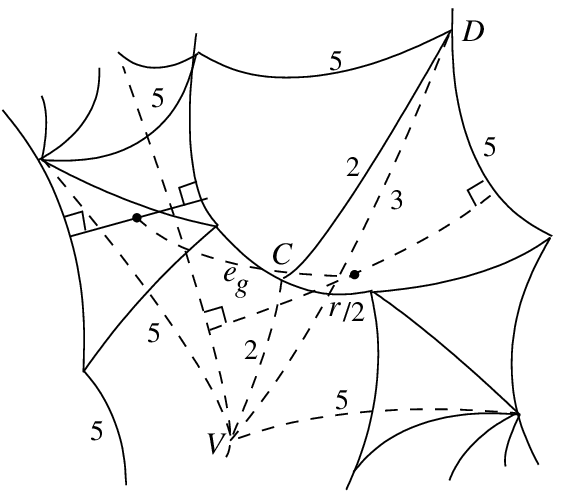}\qquad&
\qquad\includegraphics[width=4 cm]{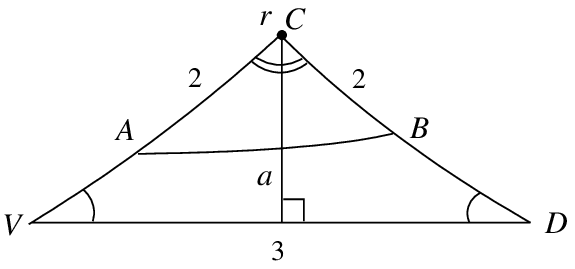}\\
(a)\qquad &\qquad (b)\\
\end{tabular}
\caption{}\label{pti_19}
\end{figure}

\medskip
In {\bf Case (19)}, $n=5$, $l=r/2$, $r\geq 3$, $(r,2)=1$,
$\Delta=(5,r,5)$. It is not
difficult to construct the decomposition of $\widetilde\P$ into
tetrahedra $T=[2,2,3;5,5,r]$ of infinite volume. The axis $e_1$
coincides with an edge of $T$. We determine the position of $e_g$
(Figure~\ref{pti_19}a).
Consider the face $F$ of $\widetilde\P$ lying in $\omega$. $F$
is split into two parts by the line $l=\eta\cap\omega$, where $\eta$
is a reflection plane passing through $e_1$ orthogonally to the
bisector $\zeta$. Since the (acute) angle that $l$ makes with $f$
is greater than the (acute) angle that $l$ makes with
$\omega\cap\alpha'$, $e_g$ intersects $\eta$ in an interior point
of $VCD$. Moreover, the angle between $e_g$ and $\eta$ is less
than $\pi/6$.

Then there is an elliptic element $h\in \Gamma^*$ of order $q>3$
with an axis lying in $\eta$ and disjoint from $VD$. Let $h$
intersect the sides $VC$ and $CD$ of $\triangle VCD$ in points $A$ and
$B$, respectively (Figure~\ref{pti_19}b). Clearly, $h$ cannot
intersect both the sides at right angles. So, $St_{\Gamma^*}(A)$
and $St_{\Gamma^*}(B)$ cannot be both $(2,2,q)$. Therefore,
$q=4$ or $5$. Note that $VCD$ is an isosceles triangle with
$$
\angle V=\angle D= \arccos ((1+\sqrt{5})/2)\quad {\rm and} \quad
\angle C=\pi/r.
$$

Suppose $r=3$. For $q=5$ calculate the altitude $a$ from $C$
in $VCD$: $\cosh a=4\cos(\pi/5)/\sqrt{3}$. But then
$a/2<\rho_{min}(3,5)$ and thus, the group is not discrete. Suppose
$q=4$. The axis $h$ of order $4$
intersects the axes of order $2$ at $A$ and $B$.
Since in a discrete group
axes of orders $2$ and $4$ intersect at $\pi/4$ or $\pi/2$ and
$ABC$ is a hyperbolic triangle, we have $\angle A=\angle B=\pi/4$
in $\triangle ABC$. Then the distance between $AB$ and $VD$ is
less than $\rho_{min}(3,4)$. Thus, $\Gamma^*$ is not discrete.

As for $r\geq 5$, it turns out that
$a-\rho_{min}(q,r)<\rho_{min}(q,3)$ for $q=4,5$. Thus, the group
is not discrete.

We have considered all the cases and so Theorem~A is proved.
\qed

\section{On the minimal volume hyperbolic $3$-orbifold}

Let $G_T$ be the group generated by reflections in the faces of the
hyperbolic tetrahedron $T=[2,2,3;2,5,3]$ and  let $\Delta_T$ be its
orientation preserving index~2 subgroup. It is well-known that
$\Delta_T$ can be generated by two elements $a$ and $b$
of order $3$ whose axes are the mutually orthogonal edges of $T$
with dihedral angles $\pi/3$ (see eg. \cite{K2}).
The ${\Bbb Z}_2$-extension
$\Gamma_{353}=\left<\Delta_T,e\right>$, where $e$ is a half-turn whose axis
is orthogonal to the opposite edges of $T$ with dihedral angles
$\pi/2$ and $\pi/5$, is discrete by~\cite{DM}. $\Gamma_{353}$ has minimal
co-volume among all Kleinian groups containing elements
of order $p\geq 4$ and among all groups with a tetrahedral subgroup~\cite{GM1}.
It is also known that $\Gamma_{353}$
has minimal co-volume in the class of arithmetic Kleinian
groups~\cite{CF}.

The group $\Gamma_{353}=\left<a,b,e\right>$ is generated by $a$ and
$e$ because $b=eae$.
For this choice of generators, the parameter $\gamma(a,e)$ is not
real, since the axes of $a$ and $e$ are not mutually orthogonal.
(For geometric meaning of the parameters
see~\cite[Theorems~1--3]{KK}.) However, the following statement is true.

\medskip
\noindent
{\bf Corollary 1}
{\it
The group $\Gamma_{353}$ is an ${\cal RP}$ group.
}

\medskip
\noindent
{\it Proof.}
The idea is to find another generating pair for $\Gamma_{353}$.
We will show that $\Gamma_{353}$ is isomorphic to one of the
discrete $\cal RP$ groups that appear in the proof of
Theorem~A.
Consider, for example, the group from Item ({\it ii}) of
Theorem~A with $n=5$, $m=2$, and $l=3/2$
(see  Case~(11) in Section~5 and Figure~\ref{ptf_discr}h).
Recall that $\widetilde\Gamma=\langle f,g,e\rangle$ is the orientation
preserving subgroup of $\Gamma^*=\langle G_T,e_g\rangle$
and so
$\widetilde\Gamma$ is isomorphic to $\Gamma_{353}$.

It remains to prove that
$\widetilde\Gamma=\Gamma$. All we need is to show that
$\Gamma_{353}\cong\widetilde\Gamma$ is generated by $f$ and $g$,
i.e. that $e=W(f,g)$.

We have
$$
h_1^2=R'_\alpha R_\alpha=gfg^{-1}f
\quad{\rm and}\quad
h_2^2=R''_\alpha R'_\alpha=f^2g^{-1}f^{-1}gf^2gf^{-1}g^{-1}.
$$
Moreover, in our case $h_1^4=1$ and $h_2^3=1$.
Let $x=R_\zeta R_\alpha$ and $y=R_\zeta R_\omega$, where
$\zeta$ is the bisector of the dihedral angle of
$\widetilde{\cal P}$ made by $\alpha'$ and $\alpha''$ .
Then using $(gfg^{-1}f)^2=1$ and $f^5=1$, we have
\begin{eqnarray*}\label{ufg}
x&=&R_\zeta R_\alpha=(R_\zeta R'_\alpha)(R'_\alpha R_\alpha)=
h_2h_1^2=gfg^{-1}f^3g^{-1}fgf^3gfg^{-1}f\\
&=&gfg^{-1}f^3g^{-1}fgf^2gf^{-1}g^{-1},\\
y&=&R_\zeta R_\omega=(R_\zeta R_\alpha)(R_\alpha R_\omega)=
xf^{-1}=gfg^{-1}f^3g^{-1}fgf^2gf^{-1}g^{-1}f^{-1}.
\end{eqnarray*}
The element $z=yx$ is an element of order $3$ and
$z=R_\alpha R_\eta$, where $\eta$ is the plane passing
through $e_1$ orthogonally to $\zeta$. So,
\begin{eqnarray*}
e_1&=&R_\zeta R_\eta=(R_\zeta R_\alpha)(R_\alpha R_\eta)=
xz\\
&=&gfg^{-1}f^3g^{-1}f^2gf^3gf^2g^{-1}f^3g^{-1}fgf^2gf^{-1}g^{-1}.
\end{eqnarray*}
Then since $e_1=f^2e$, we get
$$
e=f^3gfg^{-1}f^3g^{-1}f^2gf^3gf^2g^{-1}f^3g^{-1}fgf^2gf^{-1}g^{-1}.
$$
This completes the proof.
\qed

\section{Parameters}

Theorem~A can be reformulated in terms of the parameters
$\beta=\tr^2f-4$, $\beta'=\tr^2g-4$, and $\gamma=\tr[f,g]-2$.

\medskip
\noindent
{\bf Theorem B. }{\it
Let $f,g\in \rm PSL(2,{\Bbb C})$,\
$\beta=-4\sin^2(\pi/n)$,
\ $n\geq 3$, $(n,2)=1$, \ $\beta'>0$, and
\ $0<\gamma<-\beta\beta'/4$.
Then $\Gamma=\langle f,g \rangle$ is discrete if and only if
$(\beta,\beta',\gamma)$ is one of the triples listed in
rows 21--41 of Table~\ref{table_param} in Appendix.
}

\medskip
\noindent
{\it Proof.}
To prove the theorem it is sufficient to calculate values of the
parameters $(\beta,\beta',\gamma)$ for all discrete groups described
in Theorem~A. Since we know fundamental polyhedra for the discrete
groups, it can be done by straightforward calculation.

Since $f$ is a primitive elliptic element of order $n$,
$\beta=\tr^2f-4=-4\sin^2(\pi/n)$.

Compute $\gamma=\tr[f,g]-2$. Using the identity~(\ref{fprime}) from
Section~3 and the fact that $f=R_{\alpha^*}R_\omega$,
where $\alpha^*=R_\omega(\alpha)$, we have
$$
fgf^{-1}g^{-1}=f(f')^{-1}=
(R_{\alpha^*}R_\omega)(R_\omega R'_\alpha)=R_{\alpha^*}R'_\alpha.
$$
Note that $\alpha^*$ passes through
$f$ and makes an angle of $\pi/n$ with $\omega$. Moreover, $\alpha^*$
and $\alpha'$ are disjoint and $\delta$ is orthogonal to both
$\alpha^*$ and $\alpha'$.
Therefore, $[f,g]=fgf^{-1}g^{-1}$ is a hyperbolic element with the
axis lying in $\delta$ and the translation length $2d$, where $d$ is
the distance between $\alpha^*$ and $\alpha'$.

It is known \cite[Theorem~1 or Table~1]{KK} that for an elliptic
and hyperbolic generators with intersecting axes,
$0<\gamma<-\beta\beta'/4$. Hence, we get
$$\tr[f,g]=+2\cosh d \qquad{\rm and}\qquad \gamma=2(\cosh d-1).$$

Finally, we compute $\beta'=\tr^2g-4=4\sinh^2T$, where $T$ is the
distance between $e$ and $e_g$ that can be measured in $\omega$.

There are three essentially different fundamental polyhedra for the
discrete groups from Theorem~A: $\widetilde\P$, $\P_1$, and $\P_2$.
All other fundamental polyhedra are just smaller polyhedra in
decompositions of these three. We describe a procedure for
calculating the parameters in case of $\widetilde\P$, for other polyhedra
the procedure is similar.

\begin{figure}[htbp]
\centering
\includegraphics[width=4.5 cm]{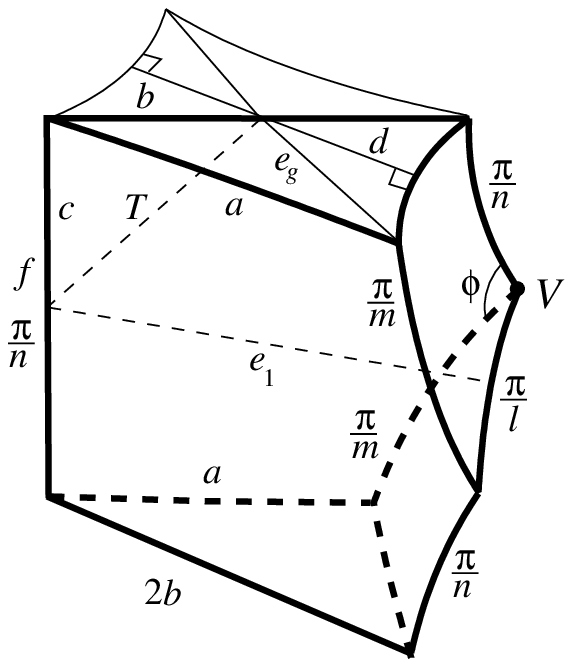}
\caption{}\label{polyhedron_par}
\end{figure}

To draw the non-compact polyhedron $\widetilde\P$  from Section~5
we add planes $\delta$ and $e_1(\delta)$, see
Figure~\ref{polyhedron_par},
where
thin lines show lines in $\delta$ that are not edges of
$\widetilde\P$, $b$ is the distance between $f$ and $e_g$, and $c$ is
the distance between $e_1$ and $\delta$.
If $m<\infty$ then from the plane $\delta$ we have
$$
\cosh d=-\cos(2\pi/n)\cos(\pi/m)+\sin(2\pi/n)\sin(\pi/m) \cosh a,
$$
where
$$
\cosh a=\frac{\cos(\pi/n)+\cos(\pi/n)\cos(\pi/m)}
{\sin(\pi/n)\sin(\pi/m)}.
$$
Analogous calculations can be done for $m=\infty$ and
$m=\overline\infty$. We obtain
$$
\gamma=\left\{
\begin{array}{lll}
2(\cos(\pi/m)+\cos(2\pi/n)) & {\rm if} & m<\infty,\\
2(1+\cos(2\pi/n)) & {\rm if} & m=\infty,\\
2(\cosh p+\cos(2\pi/n)) & {\rm if} & m=\overline\infty,\\
\end{array}
\right.
$$
where
$p$ is the distance between $\alpha$ and $\alpha'=e_g(\alpha)$ if they
are disjoint.

Further,
$$
\cosh T=\cosh b\cosh c, \quad{\rm where}\quad
\cosh b=\frac{\cos(\pi/(2m))}{\sin(\pi/n)}.
$$
Find $\cosh c$.
Suppose $V$ is proper. The face of $\widetilde\P$ lying in $\omega$ is a
pentagon with four right angles, so
$c$ is given by
$$
\cosh 2c=\frac{\cos\phi+\cosh a\cosh 2b}
{\sinh a\sinh 2b}, \quad
\cos\phi=\frac{\cos(\pi/l)+\cos(\pi/n)\cos(\pi/m)}
{\sin(\pi/n)\sin(\pi/m)}.
$$
formulas are similar for
ideal or imaginary vertex $V$ and for
$l\in\{\infty,\overline\infty\}$.

Combining the formulas above and simplifying them, we have
$$\beta'=\left\{
\begin{array}{lll}
\displaystyle\frac{2\cos(\pi/l)}{\gamma}-
 \frac{\sqrt{\beta+4}(\beta+(\gamma-\beta)^2)}{\gamma\beta}
 -\frac{2\gamma}{\beta}-2 & {\rm if} & m<\infty,\\
\displaystyle\frac{2}{\gamma}-
 \frac{\sqrt{\beta+4}(\beta+(\gamma-\beta)^2)}{\gamma\beta}
 -\frac{2\gamma}{\beta}-2
& {\rm if} & m=\infty,\\
\displaystyle\frac{2\cosh t}{\gamma}-
 \frac{\sqrt{\beta+4}(\beta+(\gamma-\beta)^2)}{\gamma\beta}
 -\frac{2\gamma}{\beta}-2
& {\rm if} & m=\overline\infty,\\
\end{array}
\right.
$$
where $t$ is the distance between $\alpha'$ and $\alpha''$.
These results lead to rows (21)--(24) in Table~\ref{table_param}.

Similarly, considering the polyhedron $\P_1$ from Section~4 we
calculate the parameters for the groups from item~(2)({\it v}) of
Theorem~A to get rows (25)--(26). Rows (27)--(28) come from
Theorem~A(2)({\it vii}) (see Section~4 and
Figure~\ref{p1_hyp}b for $\P_2$).

If a fundamental polyhedron appears as a result of decomposition of
$\widetilde\P$, $\P_1$, or $\P_2$,
then it suffices to substitute corresponding values of dihedral angles
or distances into general formulas.
So, rows (29), (32), and (35)--(40) correspond to Theorem~A(2)({\it ii});
rows (33) and (34) correspond to Theorem~A(2)({\it iii});
rows (30) and (41) correspond to Theorem~A(2)({\it iv});
the row (31) corresponds to Theorem~A(2)({\it vi}).
\qed

\section*{Appendix}

For simplicity, in the statement of Theorem~B and in Table~\ref{table_param}
all elliptic generators are assumed to be primitive.
If one (or both) generator(s) of an $\cal RP$ group is non-primitive
elliptic, Table~\ref{table_param} and Theorem~B still can be used to verify
discreteness of the group, but first we must replace
the triple $(\beta,\beta',\gamma)$,
where $\beta=-4\sin^2(q\pi/n)$, $(q,n)=1$, $1<q<n/2$,
with a new triple $(\widetilde\beta,\beta',\widetilde\gamma)$,
where $\widetilde\beta=-4\sin^2(\pi/n)$ and
$\widetilde\gamma=(\widetilde\beta/\beta)\gamma$.
The new triple corresponds to the same group by
Gehring and Martin~\cite{GM7} (cf. \cite[Remark~2, p.~262]{KK}).

\medskip
\noindent
{\bf Remark A.1.}
A part of Table~\ref{table_param} first appeared in~\cite{KK}, but unfortunately,
there was a misprint. Here we correct it.
\medskip

\begin{longtable}{|c|c|c|c|}
\caption{All truly spatial discrete ${\cal RP}$ groups
whose generators have real traces.
Here all numbers $n$, $m$, $p$ are positive integers}\label{table_param}\\
\hline
  & $\beta=\beta(f)$ & $\gamma=\gamma(f,g)$ & $\beta'=\beta(g)$ \\
\hline
\endfirsthead
\caption[]{(continued)}\\
\hline
  & $\beta=\beta(f)$ & $\gamma=\gamma(f,g)$ & $\beta'=\beta(g)$ \\
\hline
\endhead
\multicolumn{4}{|c|}%
{Both generators are elliptic, mutually orthogonal skew axes}\\
\hline
1 &  \rule[-3ex]{0ex}{7ex} $-4\sin^2\frac{\pi}n, \ n\geq 3$
  & $\displaystyle{-4\cos^2\frac{\pi}p, \atop \cos\frac{\pi}p>\sin\frac{\pi}n\sin\frac{\pi}m}$
  & $-4\sin^2\frac{\pi}m, \ m\geq 3$ \\
\hline
2 & \rule[-2ex]{0ex}{5ex} $-4\sin^2\frac{\pi}n,\ n\geq 3$
  & $(-\infty,-4]$
  & $-4\sin^2\frac{\pi}m, \ m\geq 3$\\
\hline
3 & \rule[-2.5ex]{0ex}{6ex} $\displaystyle{-4\sin^2\frac{\pi}n, \ n\geq 7, \atop (n,2)=1}$
  & $-(\beta+2)^2$
  & $\beta$\\
\hline
\multicolumn{4}{|c|}%
{$f$ is elliptic and $g$ is parabolic, the axis of $f$ lies in an invariant plane
 of $g$}\\
\hline
4 & \rule[-2ex]{0ex}{5ex} $-4\sin^2\frac{\pi}n$,\ $n\geq 3$
  & $-4\cos^2\frac{\pi}p$, \ $p\geq 3$
  & $0$ \\
\hline
5 & \rule[-2ex]{0ex}{5ex} $-4\sin^2\frac{\pi}n$,\ $n\geq 3$
  & $(-\infty,-4]$
  & $0$\\
\hline
\multicolumn{4}{|c|}%
{Each of $f$ and $g$ is either parabolic or hyperbolic,}\\
\multicolumn{4}{|c|}%
{each of the generators has an invariant plane which is orthogonal}\\
\multicolumn{4}{|c|}%
{to all invariant planes of the other generator}\\
\hline
6 & \rule[-2ex]{0ex}{5ex} $[0,+\infty)$
  & $-4\cos^2\frac{\pi}p$, \ $p\geq 3$
  & $[0,+\infty)$\\
\hline
7 & \rule[-2ex]{0ex}{5ex} $[0,+\infty)$
  & $(-\infty,-4]$
  & $[0,+\infty)$ \\
\hline
\multicolumn{4}{|c|}%
{$f$ is elliptic, $g$ is hyperbolic, disjoint axes}\\
\hline
8 & \rule[-2ex]{0ex}{5ex} $-4\sin^2\frac{\pi}n, \ n\geq 3$
  & $-4\cos^2\frac{\pi}p, \ p\geq 3$
  & $(0,+\infty)$ \\
\hline
9 & \rule[-2ex]{0ex}{5ex} $-4\sin^2\frac{\pi}n$, \ $n\geq 3$
  & $(-\infty,-4]$
  & $(0,+\infty)$ \\
\hline
10 & \rule[-3ex]{0ex}{7ex} $\displaystyle{-4\sin^2\frac{\pi}n, \ n\geq 5, \atop (n,2)=1}$
  & $-(\beta+2)^2$
  & $\displaystyle{4(\beta+4)\cos^2\frac{\pi}p-4, \atop p\geq 4}$\\
\hline
11 & \rule[-3ex]{0ex}{7ex} $\displaystyle{-4\sin^2\frac{\pi}n, \ n\geq 5, \atop (n,2)=1}$
  & $-(\beta+2)^2$
  & $[4(\beta+3),+\infty)$\\
\hline
12 & \rule[-2.5ex]{0ex}{6ex} $-3$
  & $(\sqrt{5}-3)/2$
  & $\displaystyle{2\cos^2\frac{\pi}p(7+3\sqrt{5})-4, \atop p\geq 3}$\\
\hline
13 & \rule[-1ex]{0ex}{4ex} $-3$
  & $(\sqrt{5}-3)/2$
  & $[2(5+3\sqrt{5}),+\infty)$\\
\hline
\multicolumn{4}{|c|}%
{$f$ is elliptic of even order $n$, $g$ is hyperbolic,}\\
\multicolumn{4}{|c|}%
{the axes intersect non-orthogonally; $1/n+1/m<1/2$}\\
\hline
14 & \rule[-2ex]{0ex}{6ex} $-4\sin^2\frac{\pi}n,\ n\geq 4$
   & $\displaystyle {2(\cos\frac {2\pi}m+\cos\frac {2\pi}n),
       \atop (m,2)=2}$
   & $\displaystyle \frac{4\cos^2\frac\pi p}\gamma-\frac {4\gamma}\beta$,
     $p\geq 3$\\
\hline
15 & \rule[-3ex]{0ex}{7ex} $-4\sin^2\frac{\pi}n,\ n\geq 4$
   & $\displaystyle {2(\cos\frac {2\pi}m+\cos\frac {2\pi}n) \atop (m,2)=2}$
   & $\displaystyle\left[\frac 4\gamma-\frac{4\gamma}\beta, +\infty\right)$\\
\hline
16 & \rule[-2ex]{0ex}{6ex} $-4\sin^2\frac{\pi}n,\ n\geq 4$
   & $\left[\beta+4,+\infty\right)$
   & $\displaystyle \frac{4\cos^2\frac\pi p}\gamma-\frac {4\gamma}\beta$, $p\geq 3$\\
\hline
17 & \rule[-3ex]{0ex}{7ex} $-4\sin^2\frac{\pi}n,\ n\geq 4$
   & $\left[\beta+4,+\infty\right)$
   & $\displaystyle \left[\frac 4\gamma-\frac{4\gamma}\beta, +\infty\right)$\\
\hline
18 & \rule[-2ex]{0ex}{6ex} $-4\sin^2\frac{\pi}n,,\ n\geq 4$
   & $\displaystyle {2(\cos\frac{2\pi}m+\cos\frac{2\pi}n), \atop (m,2)=1}$
   & $\displaystyle {\frac{4(\gamma-\beta)}\gamma\cos^2\frac\pi p-\frac{4\gamma}\beta, \atop
     p \geq 3}$\\
\hline
19 & \rule[-3ex]{0ex}{7ex} $-4\sin^2\frac{\pi}n,\ n\geq 4$
   & $\displaystyle {2(\cos\frac{2\pi}m+\cos\frac{2\pi}n), \atop (m,2)=1}$
   & $\displaystyle \left[\frac{4(\gamma-\beta)}\gamma-\frac{4\gamma}\beta,+\infty\right)$\\
\hline
20 & \rule[-2ex]{0ex}{6ex} $-2$
   & $\displaystyle{2\cos\frac{2\pi}m, \atop m\geq 5,\ (m,2)=1}$
   & $\gamma^2+4\gamma$\\
\hline
\multicolumn{4}{|c|}%
{$f$ is elliptic of odd order $n$, $g$ is hyperbolic, the axes intersect}\\
\multicolumn{4}{|c|}%
{non-orthogonally; $1/n+1/m<1/2$,}\\
\multicolumn{4}{|c|}%
{$\displaystyle U=-2\frac{(\gamma-\beta)^2\cos\frac{\pi}{n}+\gamma(\gamma+\beta)}{\gamma\beta}$, \
 $\displaystyle V=-2\frac{(\beta+2)^2\cos\frac{\pi}{n}}{\beta+1}-2\frac{(\beta^2+6\beta+4)}\beta$}\\
\hline
21 & \rule[-2ex]{0ex}{6ex} $-4\sin^2\frac{\pi}n, \ n\geq 3$
   & $\displaystyle {2(\cos\frac{2\pi}{m}+\cos\frac{2\pi}{n}), \atop (m,2)=2}$
   & $\displaystyle {\frac 2\gamma(\cos\frac{\pi}{p}-\cos\frac{\pi}{n})+U, \atop p\geq 2}$\\
\hline
22 & \rule[-3ex]{0ex}{7ex} $-4\sin^2\frac{\pi}n, \ n\geq 3$
   & $\displaystyle {2(\cos\frac{2\pi}{m}+\cos\frac{2\pi}{n}), \atop (m,2)=2}$
   & $\displaystyle \left[\frac{2(1-\cos\frac{\pi}{n})}{\gamma}+U,+\infty\right)$\\
\hline
23 & \rule[-2ex]{0ex}{6ex} $-4\sin^2\frac{\pi}n, \ n\geq 3$
   & $[\beta+4,+\infty)$
   & $\displaystyle {\frac 2\gamma (\cos\frac{\pi}{p}-\cos\frac{\pi}{n})+U, \atop p\geq 2}$\\
\hline
24 & \rule[-3ex]{0ex}{7ex} $-4\sin^2\frac{\pi}n, \ n\geq 3$
   & $[\beta+4,+\infty)$
   & $\displaystyle \left[\frac{2(1-\cos\frac{\pi}{n})}{\gamma}+U,+\infty\right)$\\
\hline
25 & \rule[-2ex]{0ex}{6ex} $-4\sin^2\frac{\pi}n, \ n\geq 3$
   & $\displaystyle {2(\cos\frac{2\pi}{m}+\cos\frac{2\pi}{n}), \atop (m,2)=1}$
   & $\displaystyle {\frac{2(\gamma-\beta)}{\gamma}\cos\frac{\pi}{p}+U, \ p\geq 2}$\\
\hline
26 & \rule[-3ex]{0ex}{7ex} $-4\sin^2\frac{\pi}n, \ n\geq 3$
   & $\displaystyle {2(\cos\frac{2\pi}{m}+\cos\frac{2\pi}{n}), \atop (m,2)=1}$
   & $\displaystyle \left[\frac{2(\gamma-\beta)}{\gamma}+U,+\infty\right)$\\
\hline
27 & \rule[-2.5ex]{0ex}{7ex} $-4\sin^2\frac{\pi}n, \ n\geq 7$
   & $(\beta+4)(\beta+1)$
   & $\displaystyle{\frac{2(\beta+2)^2\cos\frac{\pi}{p}}{\beta+1}+V, \ p\geq 2}$\\
\hline
28 & \rule[-3ex]{0ex}{7ex} $-4\sin^2\frac{\pi}n, \ n\geq 7$
   & $(\beta+4)(\beta+1)$
   & $\displaystyle\left[\frac{2(\beta+2)^2}{\beta+1}+V,+\infty\right)$\\
\hline
29 & \rule[-2ex]{0ex}{6ex} $\displaystyle{-4\sin^2\frac{\pi}n, \ n\geq 5, \atop (n,3)=1}$
   & $\beta+3$
   & $\displaystyle2\frac{(\beta-3)\cos\frac{\pi}n-2\beta-3}{\beta}$\\
\hline
30 & \rule[-2ex]{0ex}{6ex} $\displaystyle{-4\sin^2\frac{\pi}n, \ n\geq 5, \atop (n,3)=1}$
   & $2(\beta+3)$
   & $\displaystyle-6\frac{2\cos\frac{\pi}n+\beta+2}{\beta}$\\
\hline
31 & \rule[-2ex]{0ex}{6ex} $-3$
   & $\displaystyle{2\cos(2\pi/m)-1, \atop m\geq 7, (m,2)=1}$
   & $\displaystyle{2\frac{\gamma^2+2\gamma+2}{\gamma}}$\\
\hline
32 & \rule[-2ex]{0ex}{5ex} $-3$
   & $\displaystyle{2\cos(\pi/m)-1, \atop m\geq 4, (m,3)=1}$
   & $\gamma^2+4\gamma$\\
\hline
33 & \rule[-2ex]{0ex}{5ex} $-3$
   & $\displaystyle{2\cos(2\pi/m), \atop m\geq 7, (m,4)\leq 2}$
   & $2\gamma$\\
\hline
34 & \rule[-1ex]{0ex}{4ex} $-3$ & $(\sqrt{5}+1)/2$ & $\sqrt{5}$\\
\hline
35 & \rule[-1ex]{0ex}{4ex} $-3$ & $(\sqrt{5}-1)/2$ & $\sqrt{5}$\\
\hline
36 & \rule[-1ex]{0ex}{4ex} $-3$ & $(\sqrt{5}-1)/2$ & $\sqrt{5}-1$\\
\hline
37 & \rule[-1ex]{0ex}{4ex}$(\sqrt{5}-5)/2$ & $(\sqrt{5}-1)/2$ & $\sqrt{5}$\\
\hline
38 & \rule[-1ex]{0ex}{4ex}$(\sqrt{5}-5)/2$ & $(\sqrt{5}-1)/2$ & $(3\sqrt{5}-1)/2$\\
\hline
39 & \rule[-1ex]{0ex}{4ex}$(\sqrt{5}-5)/2$ & $(\sqrt{5}-1)/2$ & $3(\sqrt{5}+1)/2$\\
\hline
40 & \rule[-1ex]{0ex}{4ex}$(\sqrt{5}-5)/2$ & $(\sqrt{5}+1)/2$ & $3(\sqrt{5}+1)/2$\\
\hline
41 & \rule[-1ex]{0ex}{4ex}$(\sqrt{5}-5)/2$ & $\sqrt{5}+2$ & $(5\sqrt{5}+9)/2$\\
\hline
\end{longtable}

\noindent
Gettysburg College, Mathematics Department\\
300 North Washington St., Campus Box 402\\
Gettysburg, PA 17325, USA\\
{\tt yklimenk@gettysburg.edu}\\
\\
CMI, Universit\'e de Provence, 39, rue F. Joliot Curie\\
13453 Marseille cedex 13, France\\
{\tt kopteva@cmi.univ-mrs.fr}


\begin{thebibliography}{99}
%
\bibitem{B}
A.~F.~Beardon, {\it The geometry of discrete groups},
Springer-Verlag, New~York--Heidelberg--Berlin, 1983.
%
\bibitem{CF}
T.~Chinburg and E.~Friedman,
{\it The smallest arithmetic hyperbolic three-orbifold},
Invent. Math. {\bf 86} (1986), no. 3, 507--527.
%
\bibitem{DM}
D.~A.~Derevnin and A.~D.~Mednykh, {\it Discrete extensions of the Lanner groups
(Russian)}, Dokl. Akad. Nauk {\bf 361} (1998), no. 4, 439--442;
translation in: Dokl. Math. {\bf 58} (1998), no. 1, 78--80.
%
\bibitem{EP}
D.~B.~A.~Epstein and C.~Petronio, {\it An exposition of Poincar\'e's
polyhedron theorem},
L'Enseignement Math\'ematique {\bf 40} (1994), 113--170.
%
\bibitem{F2}
A.~A.~Felikson, {\it Coxeter decompositions of spherical tetrahedra},
preprint, 99--053, Bielefeld, 1999.
%
\bibitem{Fen}
W.~Fenchel,
{\it Elementary geometry in hyperbolic space},
de Gruyter Studies in Mathematics, {\bf 11}.
Walter de Gruyter \& Co., Berlin, 1989.
%
\bibitem{FR}
B.~Fine and G.~Rosenberger,  {\it Classification of generating pairs
of all two generator Fuchsian groups}, in: Groups'93 (Galway/St.Andrews,
v.1, Galway, 1993), 205--232; London Math.~Soc.~Lecture Note Ser., 211,
Cambridge Univ.~Press, Cambridge, 1995.
%
\bibitem{GGM}
F.~W.~Gehring, J.~P.~Gilman, and G.~J.~Martin,
{\it Kleinian groups with real parameters},
Commun. Contemp. Math. {\bf 3}, no.~2 (2001), 163--186.
%
\bibitem{GM7}
F.~W.~Gehring and  G.~J.~Martin,  {\it Chebyshev polynomials and
discrete groups}, Proc. of the Conf. on Complex Analysis (Tianjin,
1992), 114--125, Conf. Proc. Lecture Notes Anal., I, Internat.
Press, Cambridge, MA, 1994.
%
\bibitem{GM}
F.~W.~Gehring and  G.~J.~Martin,
{\it Commutators, collars and
the geometry of M\"obius groups}, J. Anal. Math. {\bf 63} (1994),
175--219.
%
\bibitem{GM1}
F.~W.~Gehring and  G.~J.~Martin,
{\it On the minimal volume hyperbolic $3$-orbifold},
Math. Res. Lett. {\bf 1} (1994), no. 1, 107--114.
%
\bibitem{GMM}
F.~W.~Gehring, F.~W.~Marshall, and  G.~J.~Martin,
{\it The spectrum of elliptic axial distances in Kleinian groups},
Indiana Univ. Math.~J.  {\bf 47}  (1998),  no. 1, 1--10.
%
\bibitem{Gi}
J.~Gilman, {\it Two-generator discrete subgroups of ${\rm PSL}(2,{\Bbb R})$},
Mem. Amer. Math. Soc. {\bf 117} (1995), no.~561.
%
\bibitem{GiM}
J.~Gilman and B.~Maskit,  {\it An  algorithm  for  $2$-generator
Fuchsian groups}, Mich. Math.~J. {\bf 38} (1991), no.~1, 13--32.
%
\bibitem{K2}
E.~Klimenko,  {\it Discrete  groups   in   three-dimensional
Lo\-ba\-chevsky space generated  by  two  rotations},  Siberian
Math.~J. {\bf 30} (1989), no.~1, 95--100.
%
\bibitem{K3}
E.~Klimenko, {\it Some remarks on subgroups of\/ ${\rm PSL}(2,{\Bbb C})$},
Q \& A
in General Topology {\bf 8} (1990), no.~2, 371--381.
%
\bibitem{K4}
E.~Klimenko,  {\it Some examples of discrete groups and
hyperbolic orbifolds of infinite volume},
J. Lie Theory {\bf 11} (2001), 491--503.
%
\bibitem{KK}
E.~Klimenko and N.~Kopteva, {\it Discreteness criteria for
$\cal RP$ groups},
Israel J. Math. {\bf 128} (2002), 247--265.
%
\bibitem{KK05}
E.~Klimenko and N.~Kopteva, 
{\it Discrete $\cal{RP}$ groups with a parabolic generator}, 
to appear in Siberian Math. J. in 2005.
%
\bibitem{KS}
E.~Klimenko and M.~Sakuma, {\it Two-generator discrete subgroups of
${\rm Isom}({\Bbb H}^2)$
containing orientation-reversing elements}, Geometriae Dedicata {\bf
72} (1998), 247--282.
%
\bibitem{Kn}
A.~W.~Knapp, {\it Doubly generated Fuchsian groups},  Mich. Math.
J. {\bf 15} (1968), no.~3, 289--304.
%
\bibitem{M}
B.~Maskit,  {\it Kleinian groups},  Grundlehren der  Math.~Wiss.
287, Springer-Verlag, 1988.
%
\bibitem{Ma}
J.~P.~Matelski, {\it The classification of discrete 2-generator
subgroups of ${\rm PSL}(2,{\Bbb R})$}, Israel J.~Math. {\bf 42} (1982),
no.~4, 309--317.
%
\bibitem{R}
G.~Rosenberger,  {\it All generating pairs of all two-generator
Fuchs\-ian groups}, Arch. Math. {\bf 46} (1986), no.~3, 198--204.
%
\end{thebibliography}
\end{document}